\documentclass[12pt]{article}

\usepackage{amssymb}
\usepackage{amsmath}
\usepackage{subfig}
\usepackage[ruled,vlined]{algorithm2e}
\usepackage{listings,xcolor,caption, mathtools, wrapfig}

\newcommand{\bs}[1]{\boldsymbol{#1}}

\usepackage{fullpage}

\begin{document}

\title{Parameter Estimation in Fluid Flow Models from Highly Undersampled Frequency Space Data}

\author{Miriam L\"ocke$^1$, Pim van Ooij$^2$, Crist\'obal Bertoglio$^3$\\
    $^1$\small{Bernoulli Institute, University of Groningen, The Netherlands,
    \tt{m.locke@rug.nl}}\\
    $^2$\small{Department of Radiology and Nuclear Medicine, Amsterdam University
    Medical Centers,}\\ \small{University of Amsterdam, The Netherlands,
    \tt{p.vanooij@amsterdamumc.nl}}\\
    $^3$\small{Bernoulli Institute, University of Groningen, The Netherlands,
    \tt{c.a.bertoglio@rug.nl}}\\\small{(corresponding author)}\\
}
\maketitle

\begin{abstract}
    4D Flow MRI is the state-of-the-art technique for measuring blood flow and
    provides valuable information for inverse problems in the cardiovascular system.
    However, 4D Flow MRI requires very long acquisition times, straining healthcare
    resources and inconveniencing patients. To address this, usually only a part of
    the frequency space is acquired, which necessitates further assumptions to obtain an image.

    Inverse problems based on 4D Flow MRI data have the potential to compute
    clinically relevant quantities without invasive procedures and to expand the set
    of biomarkers for more accurate diagnosis. However, reconstructing MRI
    measurements with compressed sensing techniques introduces artifacts and
    inaccuracies, which can compromise the results of inverse problems. Additionally,
    there are many different sampling patterns available, and it is often unclear
    which is preferable.

    Here, we present a parameter estimation problem that directly uses highly
    undersampled frequency space measurements. This problem is numerically solved by a
    Reduced-Order Unscented Kalman Filter (ROUKF). We show that this approach results
    in more accurate parameter estimation for boundary conditions in a synthetic
    aortic blood flow model than using measurements reconstructed with compressed sensing.

    We also compare different sampling patterns, demonstrating how the quality of
    parameter estimation depends on the choice of sampling pattern. The results show
    considerably higher accuracy than inverse problems using velocity measurements
    reconstructed via compressed sensing. Finally, we confirm these findings on real
    MRI data from a mechanical phantom.

    Keywords: 4D MRI, inverse problems, Kalman filter, finite element method
\end{abstract}

\section{Introduction}
\label{sec:introduction}
In cardiovascular modeling in general and blood flows in particular, the
personalization of spatially distributed (i.e.~3D) models is a key step in performing
predictive patient-specific simulations \cite{nolte2022review}. This requires the
estimation of relevant parameters from clinical data, which can be used in diagnostics
or in the creation of patient-specific predictive simulations such as 3D models of the
hemodynamics in the vascular system.
A common technique for acquiring measurements of blood flows is phase-contrast
Magnetic Resonance Imaging (MRI),  or PC-MRI \cite{markl2003jmri}, as it is
non-invasive and non-ionizing. However, acquiring 4D flow data with MRI requires long
acquisition times, which strains clinical resources and inconveniences patients. This
applies especially if a high spatial or temporal resolution in the data is required.

MRI applies sequences of magnetic fields and measures the precession of the quantum
spin of hydrogen atoms in a strong magnetic field. The resulting measurement equates
to the Fourier transformation of the complex magnetization, the phase of which
contains information about the velocity. As such, MRI scans in frequency space (also
called k-space).
One of the features of MRI is that it is possible to select which frequencies of the
k-space to measure. Therefore, to reduce the long acquisition times, often only a
small part of the k-space is acquired \cite{gottwald2020pseudo}.

This partial acquisition introduces artifacts in the reconstructed image that reduce
the image quality. As a result, various Compressed Sensing (CS) techniques have been
developed, which aim to reconstruct velocities from highly undersampled data while
minimizing artifacts.
CS has allowed for more important reductions in the sampling, in the sense that the
accuracy of the reconstructed velocity images is higher than with the previously
developed sampling and reconstruction strategies \cite{peper2020highly, cs2007,
ye_compressed_2019}.

Nonetheless, as we will show later in this article, using CS-reconstructed data for
parameter estimation compromises the accuracy of the estimated flow and parameters
considerably. Secondly, the potential of directly using the k-space measurements in
the inverse problem to estimate parameters in blood flow has remained unexplored,
while having the potential to avoid the introduction of reconstruction artifacts.

In this paper, we introduce a new technique for estimating blood flow parameters
directly from the undersampled k-space MRI data by using an objective function
designed for complex MRI measurements. This both avoids artifacts and improves
efficiency by skipping the reconstruction/compressed sensing step. Preliminary results
in this direction have been shown in \cite{locke_parameter_2025-1} without being
explored in-depth. This paper provides additional notes on the methodology as well as
further results for different $V_{enc}$ values and the effects per parameter, as well
as exploring the robustness of the method to using approximations of the magnitude and
background phase. Furthermore, we are also including results on real MRI phantom data.

\section{Theory}

\subsection{MRI measurement and reconstruction in a nutshell}
\label{sec:measurement}

\subsubsection{Velocity encoding}
\label{sec:mri_meas}
Let us denote by $u(\bs{x}, t)$ the component of the velocity $\bs{u}(\bs{x}, t)$ in
the direction $\bs{d}$ (fixed in space and in time), and $\bs{x}\in \mathbb{R}^3$,
$t\in \mathbb{R}$ represent the spatial and temporal location within the image space,
respectively.
The MR images usually analyzed in the clinical context are actually complex valued,
namely the so-called \textit{magnetization}
\begin{equation}
    m(\bs{x}, t) =M(\bs{x}, t) \exp{(i\phi(\bs{x}, t))} \in \mathbb{C},
\end{equation}
where the magnitude $M > 0$ usually displays the anatomy, and the phase $\phi(\bs{x},
t) \in (-\pi, \pi]$ in blood flow imagining takes the form
\begin{equation}
    \phi(\bs{x}, t) = \frac{\pi}{venc}u(\bs{x}, t) + \phi_{back}(\bs{x}, t)
\end{equation}
where $\phi_{back}$ is the background phase, and $venc\in \mathbb{R}$ is the selected
velocity encoding (which is a function inversely proportional to the strength and
duration of the velocity encoding magnetic gradients).

Since both $u$ and $\phi_{back}$ are unknown at every voxel, at least two measurements
are required. In the simplest case (which is the most commonly used in clinical
practice), one measurement with no velocity encoding gradient is acquired to obtain
$\phi_{back}$, and a second one with $venc \neq 0$, leading to two complex magnetic
measurements, so that the velocity can be simply obtained by subtracting the phases, namely,
\begin{equation}\label{eq:pcmri}
    \hat{u} = \frac{\phi - \phi_{back}}{\pi}venc .
\end{equation}
Therefore, this technique is called \textit{Phase-Contrast MRI (PC-MRI)}.

In practice, the image $m(\bs{x}, t)$ is discrete in both space and time, being
divided into voxels and time instants. A \textit{voxel} refers to a spatial unit
across which the signal is assumed to be constant to provide a single value. The
number of voxels is determined by the spatial resolution. The time instants are the
times at which a measurement is made, which are determined by the temporal resolution.
Therefore $m(\bs{x}, t)$ can be represented as a matrix in $\mathbb{C}^{N_x\times N_y
\times N_z \times N_T}$, with $N_x, N_y, N_z$ being the number of voxels in each
spatial direction and $N_T$ the number of time instants. We will refer to the total
number of voxels as $N = N_x \cdot N_y \cdot N_z$.

\subsubsection{The raw signal and image reconstruction}\label{sec:sampling}

The raw signal measured by the MRI scanner corresponds to the spatial Fourier
transform of the magnetization, which can, at each measured ``time instant", be formulated as
\begin{equation}\label{eq:freq_meas}
    \bs{Y}(\bs k) = \mathcal{F}\left[\bs{M}\odot
        \exp\left(i\left(\frac{\pi}{venc}\bs{u}_{meas} +
    \phi_{back}\right)\right)\right](\bs k) + \epsilon(\bs k)
\end{equation}
where $\odot$ denotes the Hadamard product. We now consider discrete quantities
$\bs{M} \in \mathbb{R}^N$, $\bs{u}_{meas}\in \mathbb{R}^N$  and $\mathcal{F}:
\mathbb{R}^N \to \mathbb{R}^N$ is the three-dimensional discrete Fourier transform defined as
\begin{equation}
    \mathcal{F}[\bs{X}](\bs{k}) = \sum_{n_x=0}^{N_x-1}\sum_{n_y=0}^{N_y-1}\sum_{n_z =
    0}^{N_z - 1} \bs{X}_{n_xn_yn_z} \exp\left(-i2\pi \left(\frac{k_1n_x}{N_x} +
    \frac{k_2n_y}{N_y} + \frac{k_3n_z}{N_z}\right)\right)dx
\end{equation}
and $\bs{\epsilon} (\bs k)\in \mathbb{C}^N$ is a complex Gaussian noise with a mean of
zero\cite{IRARRAZAVAL2019}. Hence the measurements are $\bs{Y}^n \in \mathbb{C}^N$ for
$n = 1, \cdots, N_T$.

In practice, for a given temporal resolution, only a number of k-space lines can be
acquired. Therefore, the MRI scanner acquires different lines in k-space during each
cardiac cycle. Therefore a ``time instant" $n$ actually contains frequencies measured
during different cardiac cycles, under the assumption that there are no significant
differences between consecutive cardiac cycles (or if there are differences, for
instance in the heart rate, such data is rejected). There is also freedom in choosing
different k-space trajectories at different instants of the cardiac cycle, though that
is rarely done.

The noise level depends on the setup chosen in the scans. For instance, the lower the
spatial resolution (i.e., the larger the voxel size), the higher the signal-to-noise
ratio (because the signal increases, for the same noise level of the device). Also,
the lower the $venc$ the higher the sensitivity of the magnetization phase to the
velocity, and therefore one expects lower noise in the reconstructed velocities with
phase-contrast. As the phase is limited to the range $[-\pi, \pi)$, the velocity is
however limited to the range $[-venc, venc)$. If the maximum velocity exceeds the
chosen $venc$, this results in phase aliasing artifacts. As a result, a good choice of
$venc$ necessitates a compromise between noise level and presence of artifacts.

\bigskip
\noindent\textit{Undersampling and reconstruction.} If the frequency space is fully
sampled, i.e.~the values of all k-space locations are known, the velocity can be
reconstructed by applying the inverse Fourier transform to $\bs{Y}^n$ and equation
\eqref{eq:pcmri}. However, fully sampling the k-space is not practicable in
cardiovascular MRI scans since 3D scans may take of the order of hours, depending on
the spatiotemporal resolution chosen.

Therefore, to reduce the length of the scan, the frequency space is generally
undersampled. There are a large number of k-space  sampling patterns used in MRI
applications, usually in combination with a reconstruction algorithm specific to the
sampling pattern.

Some common 2D patterns include regular Cartesian sampling, radial sampling, and
pseudo-spiral sampling \cite{feng2022review}. These can be extended into 3D by
stacking 2D designs, leading to the stack-of-stars and stack-of-spirals designs
\cite{feng2022review}. The patterns can differ in whether lines are arranged in a
horizontal and vertical or in a shear pattern, in the number of radial lines or spiral
arms, the number of points per line, and whether the patterns are stacked with
equidistant or variable density in the z-direction. Additionally, with stack-of-stars
and stack-of-spirals, there is the option of rotating the pattern by a certain angle
for each level in the stack.

Another option is to use pseudo-random sampling methods, usually with a higher density
in the center, to ensure that the artifacts are incoherent and resemble noise
\cite{peper2020highly}. Here as well there are many choices for which probability
distribution to use.

Undersampled data cannot be reconstructed with an inverse Fourier transform, as it
violates the Nyquist limit. \textit{Compressed sensing} (CS) techniques are commonly
used to reconstruct undersampled data in a way that limits the presence of artifacts.
Generally this can be formulated as solving the optimization problem
\begin{equation}
    \bs{\hat{m} = }\text{argmin}_{\bs{m}}\{||\mathcal{F}_U(\bs{m}) - \bs{Y}||^2_2 +
    \lambda ||\Phi \bs{m}||_1\}
\end{equation}
where $\mathcal{F}_U$ is the Fourier transform combined with the undersampling,
$\bs{Y}$ is the measured k-space data, $\bs{\hat{m}}$ is the reconstructed image data,
$\Phi$ is a sparsifying transform, and $\lambda$ is the regularization parameter.

The sparsifying transform transports the image into a domain where the desired image
is sparse. For MRI images, several such domains exist, for example wavelet domains or
total variation in time. A good choice of $\Phi$ depends on the spatial and temporal
resolution as well as the nature of the artifacts. This makes pseudo-random masks
appealing, as they lead to incoherent noise artifacts. The regularization parameter
$\lambda$ is usually determined heuristically or empirically and is known to have an
impact on the quality of the reconstructed images.

\subsection{The k-space parameter estimation problem}
\label{sec:inverse_problem}

\subsubsection{Parameters estimation in PDEs}

To consider the underlying physics of the the application, we assume that they can be
modelled as the solution $\bs{u}: \mathbb{R}^3\times [0, T] \to \mathbb{R}^3$ of a
partial differential equation
\begin{equation}
    F\left(\bs{u}, \frac{\partial \bs{u}}{\partial x_1}, \frac{\partial
        \bs{u}}{\partial x_2}, \frac{\partial \bs{u}}{\partial x_3}, \frac{\partial
    \bs{u}}{\partial t}, \bs{\theta}\right) = 0
\end{equation}
in a domain $\Omega \subset \mathbb{R}^3$ with an initial condition $\bs{u}(\bs{x}, 0)
= \bs{u}^0$ and a set of boundary conditions which is dependent on a set of parameters
$\bs{\theta}$, which can describe boundary conditions or material parameters. Then we
can define the forward operator $
\mathcal{A}(\bs{\theta})$ which describes the solution of this PDE for the vector of
model parameters $\bs{\theta}\in \mathbb{R}^p$. As such, the forward problem generates
data according to a physical model with given model parameters.

The goal of the inverse problem overlying this PDE is to estimate a (sub)set of the
parameters $\bs{\theta}$ given measurements of $\bs{u}(\bs{x}, t)$.

\subsubsection{Measurements}

For formulating the inverse problem to estimate the parameters $\bs{\theta}$, we have
three different choices for the measurements:
\begin{enumerate}
    \item the frequency space data $\bs{Y}^n \in \mathbb{C}^N$ for $n = 1, \cdots,
        N_T$, which is complex-valued, raw data which might be undersampled. The noise
        on this data is Gaussian, independently distributed, and zero-mean.
    \item the magnetization data $\bs{M}_{meas}^n \in \mathbb{C}^N$ for $n=1, \cdots,
        N_T$, which relates to the frequency data as $\bs{M}^n_{meas} =
        \mathcal{F}^{-1}(\bs{Y}^n) = \bs{M}_{true} + \mathcal{F}^{-1}(\bs{\epsilon})$.
        This is also complex-valued data with a Gaussian noise with zero mean, as the
        inverse Fourier transform of a Gaussian is Gaussian as well. However, if the
        k-space is undersampled, this measurement has to be reconstructed with CS. In
        this case, the noise on the velocity measurements cannot assumed to be
        Gaussian and independently distributed, as there are correlations between the
        noise in different pixels. \cite{partin2022analysisreconstructionnoiseundersampled}
    \item the velocity data $\bs{V}^n \in \mathbb{R}^N$ for $n=1, \cdots, N_T$.
        $\bs{V}^n$ is calculated from the angle of the magnetization, meaning that in
        the case of undersampling the magnetization has to be reconstructed first
        using compressed sensing as well. Additionally, the noise on this measurement
        cannot be approximated well by a Gaussian distribution\cite{IRARRAZAVAL2019}.
\end{enumerate}
For all three types of measurements, each velocity direction is acquired separately,
therefore it is possible that only one or two velocity components are available rather
than all three, or that they are acquired in a direction that does not match a vector
component in the canonical basis.

\subsubsection{Objective function}

\label{sec:roukf}
In a Bayesian framework, the inverse problem for the estimation of parameters from
velocity measurements can be solved by minimizing the functional
\begin{equation}\label{eq:inv}
    \bs{\hat{\theta}} = \text{argmin}_{\bs{\theta} \in \mathbb{R}^p}
    \frac{1}{2\sigma_v^2}\sum_{n=1}^{N_T}\sum_{s=1}^N \left(\left[\bs{V}^n -
    \mathcal{H}(\bs{u}_{\bs{\theta}}^n)\right]_s\right)^2 + \frac{1}{2} ||\bs{\theta}
    - \bs{\theta}^0||_{(\bs{P}^0)^{-1}}^2
\end{equation}
where $\bs{V}^k \in \mathbb{R}^N$ are the measurements of the velocity provided by the
reconstruction of the PC-MRI acquisition and
$\mathcal{H}: [H_1(\Omega)]^3 \to \mathbb{R}^N$ is the observation operator, which is
applied to the result $\bs{u}_{\bs{\theta}}^n$ of the forward model for a set of model
parameters $\bs{\theta}\in \mathbb{R}^p$ at the time corresponding to the measurement
$n$. In this standard inverse problem for velocities, this observation operator
usually corresponds to the interpolation from the model's (usually finer) mesh to an
array of measured velocities at each image voxel. In fact, $[.]_s$ denotes the $s$-th
vector element (or voxel in the case of image data), therefore summing over all the
vector elements. $\bs{\theta}^0$ is the initial guess for the parameters with its
covariance matrix $\bs{P}^0$, both of which are given by the user but are assumed to
be known for the numerical inverse solver. $\sigma_v$ denotes the standard deviation
of the noise on the velocity measurements, which we also assume to be (approximately) known.

This inverse problem has been commonly applied to the estimation of boundary
conditions in hemodynamics. For a thorough review, see \cite{nolte2022review}.

In \cite{garay2022}, an alternative formulation of this objective function was
proposed to account for aliasing artifacts present in velocity MRI data. The
functional in this case is
\begin{eqnarray}
    \bs{\hat{\theta}} & = & \text{argmin}_{\bs{\theta} \in \mathbb{R}^p}  \frac{1}{2
    \sigma_M^2} \sum_{n=1}^{N_T}\sum_{s=1}^N \Big(\Big(  [\Re(\bs{M}_{meas}^n)
                \notag\\ & & -  |\bs{M}_{meas}^n|\cos(\phi_{back} +
        \mathcal{H}(\bs{u}_{\bs{\theta}}^n)\frac{\pi}{venc})]_s  \Big)^2  \notag \\ &
        & +  \Big([ \Im(\bs{M}_{meas}^n)  -  |\bs{M}_{meas}^n|\sin(\phi_{back} +
    \mathcal{H}(\bs{u}^n_{\bs{\theta}})\frac{\pi}{venc} ]_s) \Big)^2\Big) \notag\\ & &
    + \frac{1}{2} ||\bs{\theta} - \bs{\theta}^0||_{(\bs{P}^0)^{-1}}^2  \notag \\
    &  = & \text{argmin}_{\bs{\theta} \in \mathbb{R}^p}  \sum_{n=1}^{N_T}\sum_{s=1}^N
    \frac{|\bs{M}_{meas}^n|_s^2}{\sigma_M^2} \Big(  1  -
        \cos(\frac{\pi}{venc}([\bs{V}^n - \mathcal{H}(\bs{u}^n_{\bs{\theta}})))]_s \Big
    )\notag\\ & &  + \frac{1}{2} ||\bs{\theta} - \bs{\theta}^0||_{(\bs{P}^0)^{-1}}^2
    \label{eq:jeremias}
\end{eqnarray}
where $\bs{M}^n_{meas}$ are measurements of the complex magnetization and $\sigma_M$
is the standard deviation of the noise for each of the (complex magnetization)
real/imaginary measurements. This formulation assumes that the magnitude
$|\bs{M}^n_{meas}|$ of the magnetization is known in order to formulate the problem in
terms of the measured and observed velocities.

Still, the aforementioned cost functions do not account for artifacts originating from
the frequency undersampling of the data. Therefore, as we will observe in the
numerical examples, the error in the data and hence in the estimated parameters in the
inverse problem drastically grows when increasing the undersampling. Here, we propose
to solve this issue by formulating the parameter estimation problem as it is done in
CS through a data fidelity term in terms of the k-space, therefore leading to the
minimization problem:

\begin{eqnarray}
    \bs{\hat{\theta}} & = & \text{argmin}_{\bs{\theta} \in \mathbb{R}^p}
    \frac{1}{2\sigma_y^2}\sum_{n=1}^{N_T}\sum_{s=1}^N ([\Re(\bs{Y}^n -
        \mathcal{H}_{\mathcal{F}}(\bs{u}_{\bs{\theta}}^n)]_s)^2 + ([\Im(\bs{Y}^n -
        \mathcal{H}_{\mathcal{F}}(\bs{u}_{\bs{\theta}}^n))]_s)^2 \notag\\ & &+
        \frac{1}{2} ||\bs{\theta} - \bs{\theta}^0||_{(\bs{P}^0)^{-1}}^2
        \label{eq:freqJ}
    \end{eqnarray}
    with the observation operator $\mathcal{H}_{\mathcal{F}}$ being defined as
    \begin{equation}
        \mathcal{H_{\mathcal{F}}}(u) = \mathcal{F}\left(\bs{M}\odot
        e^{i\frac{\pi}{venc}\bs{u} + \phi_{back}}\right) \odot \bs{S}
    \end{equation}
    where $\bs{S} \in \mathbb{R}^N$ is the sampling mask, i.e.~with the entries
    corresponding to the sampled voxels set to $1$, and the others to $0$.
    The $\bs{Y}^n$ are measurements in frequency space and $\sigma_y$ is the standard
    deviation of the noise on the frequency-space measurements, which is assumed to be
    Gaussian and zero-mean. The use of this observation operator requires that the
    magnitude of the magnetization as well as the background phase, or an
    approximation thereof, is known.

    \subsubsection{Parameter estimation with Kalman filtering}
    \label{sec:implementation}

    There are various approaches to solving the optimization problems described above.
    Adjoint-based variational data assimilation approaches fit the entire ensemble of
    measurements to the ensemble of the observations at the matching time steps.
    However, this poses significant storage requirements, as the entire trajectory
    needs to be stored. On the other hand, sequential data assimilation assimilates
    the information from each new measurement at the point that measurement occurs
    during one forward pass of the model. This leads to a sequential improvement of
    the estimation over the duration of the inverse problem as more measurements are considered.

    We have chosen to use the Reduced Order Unscented Kalman Filter
    (ROUKF)\cite{moireau-chapelle-11} in order to solve the optimization problem
    \eqref{eq:freqJ}, as it is computationally tractable and has been successfully
    used in blood-flow problems and other time-resolved problems due to its
    recursivity. ROUKF has been successfully employed for parameter estimation from
    (already reconstructed) MR images in this context before, alas not with the raw
    frequency space data \cite{arthurs2020flexible, imperiale_sequential_2021}.

    ROUKF is a sequential parameter estimator which corrects the posterior
    distribution of the state and parameters in each time step with the available
    measurement. In order to do so, ROUKF generates a set number of particles sampled
    from the prior distribution, and propagates each of them through the forward
    problem. The result of the propagation for each of the particles' state is used to
    calculate the correction of the mean and the covariance of the parameters at each
    time measurements are available.

    As a result, the filter relies on the measurement error or \textit{innovation},
    which is proportional to the derivative of the data fidelity term of the objective
    function with respect to $u$. For the objective function \eqref{eq:inv}, this
    leads to the innovation
    \begin{equation}\label{eq:standard_innov}
        \bs{\Gamma}^n = \bs{V}^n - \mathcal{H}(\bs{u}^n)
    \end{equation}

    Similarly, for the objective function in \eqref{eq:jeremias}, the innovation
    results from the derivative of the cost function with respect to the state, namely
    \begin{equation*}
        \bs{\Gamma}^n = \frac{1}{\sqrt{2}}|\bs{M}(t^n)| \sin\left(\frac{\pi}{venc}
        \cdot \left( \bs{V}^n - \mathcal{H}(\bs{u}^n)\right)\right)
    \end{equation*}
    where the factor $\frac{\sqrt{2}venc}{\pi}|\bs{M}(t^n)|^{-1}$ is added to ensure
    the equivalence of this innovation with Equation \eqref{eq:standard_innov} in the
    case of high $venc$. For more details, see \cite{garay2022}.

    As done for the first cost function, assuming a Gaussian distribution in the noise
    of the  complex magnetization measurements, the innovation for the k-space based
    objective function \eqref{eq:freqJ} is defined by:
    \begin{equation}
        \bs{\Gamma}_n =
        \begin{bmatrix}
            \Re(\bs{Y}^n) - \Re(\mathcal{H}_{\mathcal{F}}(\bs{u}^n))\\
            \Im(\bs{Y}^n) - \Im(\mathcal{H}_{\mathcal{F}}(\bs{u}^n))
        \end{bmatrix} .
    \end{equation}

    Moreover, we remark that we are adapting the Kalman filter as described in
    \cite{moireau-chapelle-11}. Instead of simplex sigma points (which involve $p +1 $
    solutions of the forward problem), we are using canonical sigma points, which
    involve $2 p$ solutions of the forward problem but have shown in our simulations
    better performance than their simplex counterparts. In particular, the results
    using both types of sigma points deviate more the larger the number of parameters
    to be estimated. Also, the results with simplex points change depending on the
    enumeration of the parameters during the estimation process, which is not the case
    with canonical points.

    Additionally, to ensure positivity of the parameters, we are reparameterizing such
    that $\bs{\theta} = \bs{\theta}^0 2^\nu$ with $\bs{\theta}^0$ the initial guess
    for the parameters and the filter being applied to $\nu$, as is also done in
    \cite{garay2022}. The initial value of $\nu$ is set to $\bs{0}$.

    \section{Methods}

    \subsection{Synthetic data}

    \subsubsection{Forward problem setup}

\begin{figure}[hbtp!]
        \centering
        \includegraphics[width=0.37\textwidth]{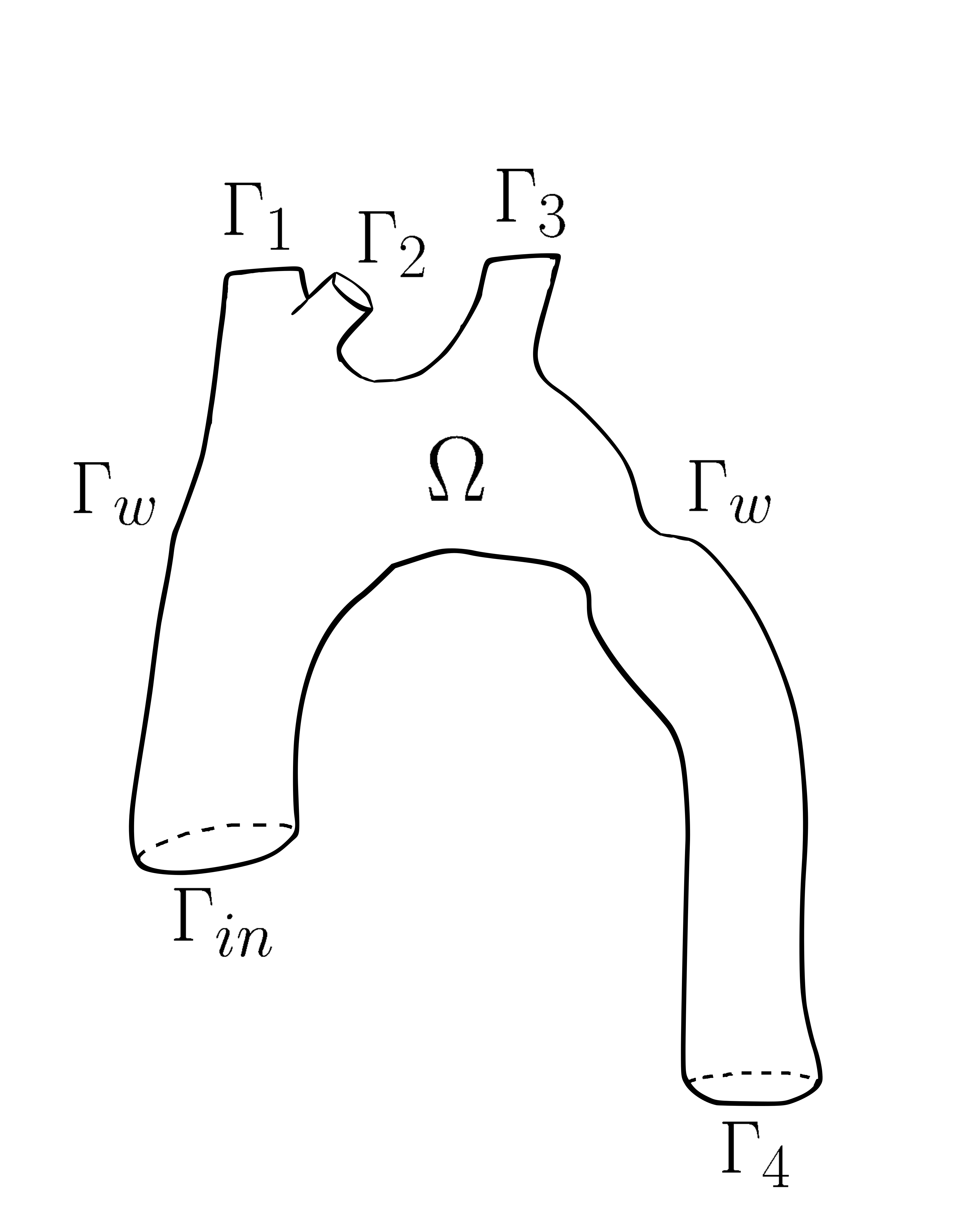}
        \caption{3D aortic model geometry}
        \label{fig:aorta}
    \end{figure}

    We use the same model as is used in \cite{garay2022}.
    We consider a geometry of the lumen of the ascending and descending aorta
    including the outlets of the brachycephalic artery, left common carotid artery,
    and left subclavian artery, as depicted in Figure \ref{fig:aorta}. This geometry
    serves as the domain for the forward model. The geometry was discretized with
    unstructured trapezoidal elements with a total of 20,916 points.

    The boundary of the geometry consists of six different boundaries: $\Gamma_{in}$
    being the inlet boundary in the ascending aorta, $\Gamma_{w}$ the arterial wall,
    and the remaining boundaries $\Gamma_l$ for $l = 1,
    \cdots, 4$ representing the outlets.

    We model the blood flow in this domain with the incompressible Navier-Stokes
    equations for the velocity $\bs{u}(\bs{x}, t)\in \mathbb{R}^3$ and the pressure
    $p(\bs{x},t)\in \mathbb{R}$:
    \begin{equation}
        \begin{cases}
            \rho \frac{\partial \bs{u}}{\partial t} + \rho (\bs{u} \cdot \nabla)\bs{u}
            - \mu \Delta \bs{u} + \nabla p = 0 \text{ in } \Omega\\
            \nabla \cdot \bs{u} = 0 \text{ in } \Omega\\
            \bs{u} = \bs{u}_{in} \text{ on } \Gamma_{in}\\
            \bs{u} = \bs{0} \text{ on } \Gamma_w\\
            \mu \frac{\partial \bs{u}}{\partial \bs{n}} - p\bs{n} = -P_l(t)\bs{n}
            \text{ on } \Gamma_l
        \end{cases}
        \label{eq:forward}
    \end{equation}
    with $\rho, \mu$ the density and dynamic viscosity of the fluid and $P_l(t)$ being
    given by a Windkessel boundary condition defined by:
    \begin{equation}
        \begin{cases}
            P_l = R_{p,l}Q_l + \pi_l\\
            Q_l = \int_{\Gamma_l} \bs{u} \cdot \bs{n} dx\\
            C_{d, l} \frac{d\pi_l}{dt} + \frac{\pi_l}{R_{d, l}} = Q_l
        \end{cases}
    \end{equation}
    This boundary condition models the effects of the remaining vascular system on the
    outlet via the proximal and distal resistances $R_p$, $R_d$ of the vasculature and
    the distal compliance $C_d$ of the vessels.

    The inflow $\bs{u}_{in}$ is defined as
    \begin{equation*}
        \bs{u}_{in} = -U f(t) \bs{n}
    \end{equation*}
    where $U$ is a constant amplitude and
    \begin{equation*}
        f(t) =
        \begin{cases}
            \sin(\frac{\pi t}{T}) \text{ if } t \leq T\\
            \frac{\pi}{T}(t-T)\exp^{-k(t-T)} \text{ if } T_c > t > T
        \end{cases}
    \end{equation*}
    with $T_c = 0.8$ and $T = 0.36$.

    The physical parameters are set as seen in Table \ref{tab:Parameters}.
    \begin{table}[!hbtp]
        \centering
        \begin{tabular}{|l|l|}
            \hline
            Parameter & Value \\
            \hline
            $\rho \ (gr \cdot cm^3)$  & $1.2$\\
            $\mu \ (P)$  &  $0.035$\\
            $U \ (cm \cdot s^{-1})$  &  $75$ \\
            $T_c \ (s)$ &  $0.80$\\
            $T \ (s)$ & $0.36$\\
            $\kappa \ (s^{-1})$ & $70$ \\
            \hline
        \end{tabular}
        \quad
        \begin{tabular}{|l|l|l|l|l|}
            \hline
            & $ \Gamma_{1}$   & $\Gamma_{2}$ & $\Gamma_{3}$ & $\Gamma_{4}$ \\
            \hline
            $R_p \ (dyn \cdot s \cdot cm^{-5})$   & $480$ & $520$  & $520$ & $200$\\
            $R_d \ (dyn \cdot s \cdot cm^{-5})$   & $7200$ & $11520$ & $11520$ & $4800$\\
            $C \ (dyn^{-1} \cdot cm^5 )$   & $4\cdot 10^{-4}$ & $3\cdot 10^{-4}$ &
            $3\cdot 10^{-4}$ & $4\cdot 10^{-4}$\\
            \hline
        \end{tabular}
        \caption{Physical parameters and numerical values of the three-element
            Windkessel parameters for every outlet.
        } \label{tab:Parameters}
    \end{table}
    The forward problem is solved using an in-house finite elements solver, with a
    semi-implicit 3D-0D coupling scheme as in \cite{garay2022} and using $P1$ elements
    for both the velocity and the pressure. The full algorithm is detailed in the appendix.

    \subsubsection{Synthetic measurements}

    The forward solution is generated with a time step of $dt = 1ms$ and undersampled
    in time to $dt_{meas} = 15ms$, leading to a total of 56 measurements. From the
    solution of the forward problem, we simulate a PC-MRI acquisition by subsampling
    into a rectangular measurement mesh with a resolution of $[2mm, 2mm, 2mm]$ and
    then applying the process described in Section \ref{sec:measurement} with a $venc$
    of double the maximal velocity. The magnitude is modelled as
    \begin{equation}
        M(\bs{x}) =
        \begin{cases}
            1.0 \text{ if } \bs{x} \text{ is in the lumen of the vessel}\\
            0.5 \text{ otherwise.}
        \end{cases}
    \end{equation}
    and the background phase was set to an arbitrary constant value of $\phi_{back} =
    7.5\cdot 10^{-2}rad$.
    Finally a complex Gaussian noise $\bs{\epsilon} \in \mathbb{C}^N$ is added with a
    signal-to-noise ratio (SNR) of 15. Fifty independent realizations of the noise
    were generated.

\begin{figure}[!htbp]
        \centering
        \subfloat[Gaussian mask, $R = 8$, $\lambda = 0.01$]{
            \includegraphics[trim=180 40 60 30, clip, width=0.3\textwidth]{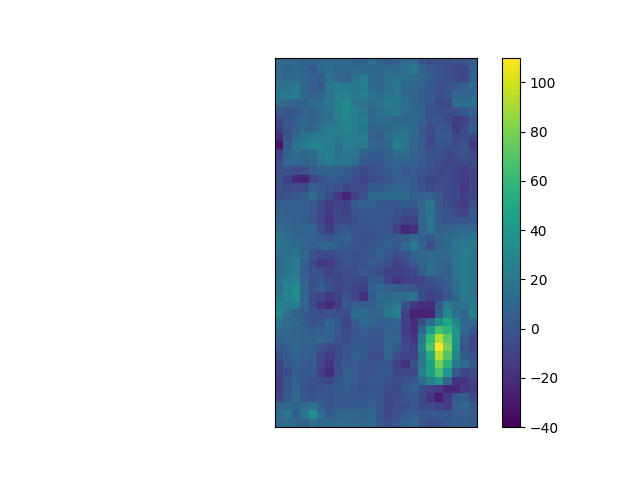}
        }
        \subfloat[Gaussian mask, $R = 16$, $\lambda = 0.01$]{
            \includegraphics[trim=180 40 60 30, clip, width=0.3\textwidth]{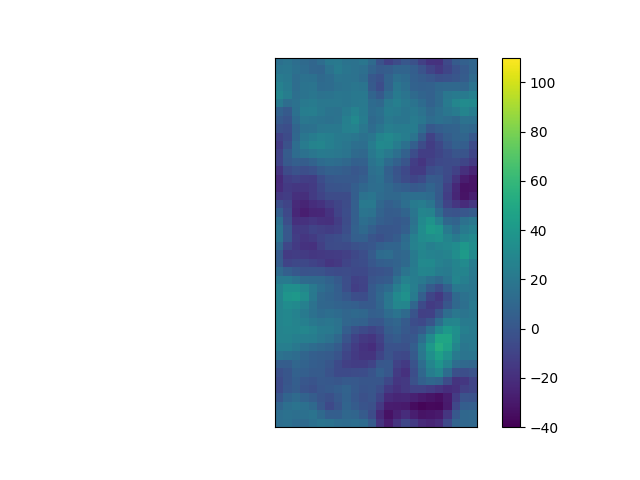}
        }
        \subfloat[Gaussian mask, $R = 32$, $\lambda = 1.0$]{
            \includegraphics[trim=180 40 60 30, clip, width=0.3\textwidth]{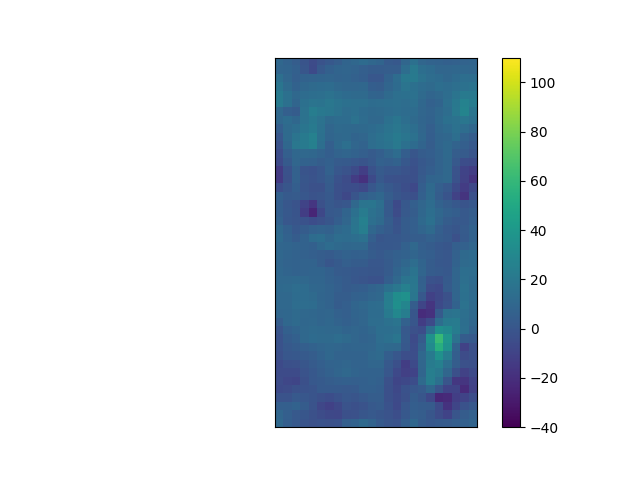}
        }  \\
        \subfloat[Spiral mask, $R = 8$, $\lambda = 0.1$]{
            \includegraphics[trim=180 40 60 30, clip, width=0.3\textwidth]{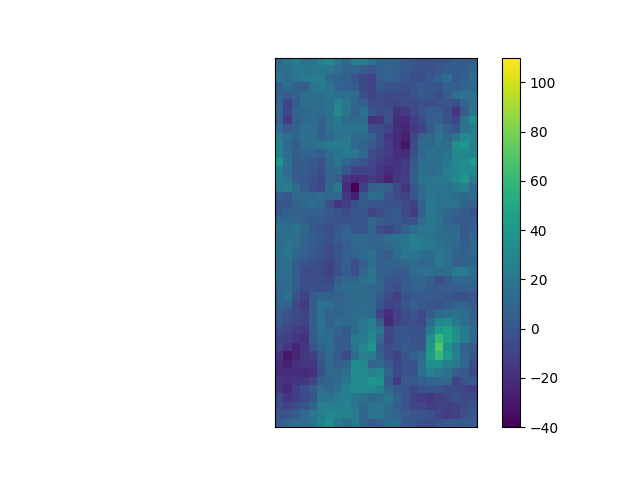}
        }
        \subfloat[Spiral mask, $R = 16$, $\lambda = 1.0$]{
            \includegraphics[trim=180 40 60 30, clip, width=0.3\textwidth]{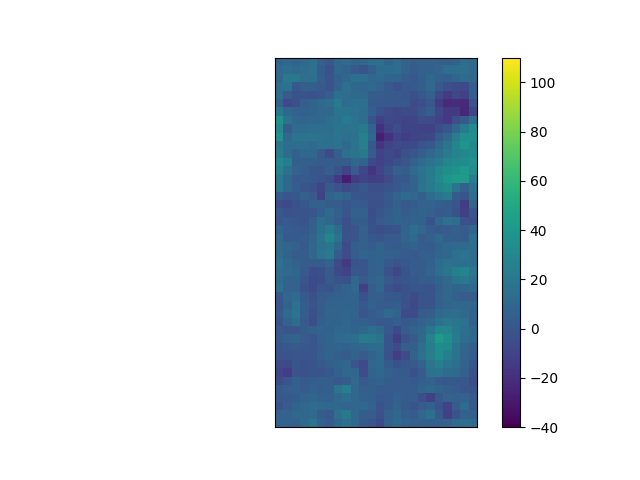}
        }
        \subfloat[Spiral mask, $R = 32$, $\lambda = 10.0$]{
            \includegraphics[trim=180 40 60 30, clip, width=0.3\textwidth]{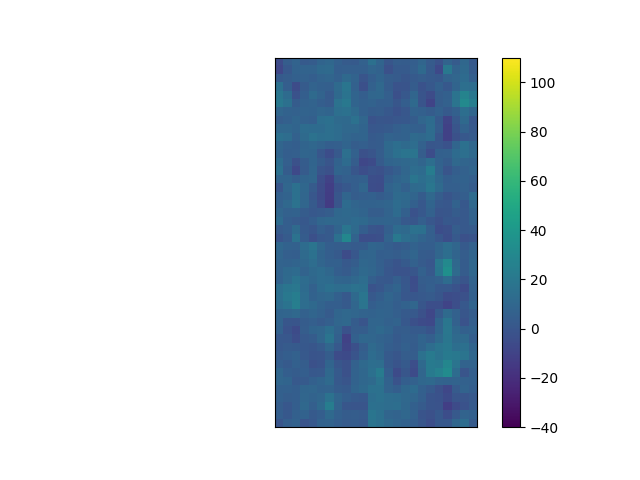}
        } \\
        \caption{Examples of velocities reconstructed with BART using different masks
            and acceleration factors. Depending on the sampling mask, different kinds of
        artifacts appear in the reconstructed velocity.}
        \label{fig:velocity_examples}
    \end{figure}

    For comparison, we reconstructed velocity measurements from these synthetic
    measurements using the Berkeley Advanced Reconstruction Toolbox
    (BART)\cite{bart2015software}. BART  is a command-line-based software that
    provides a flexible framework of compressed sensing methods, as well as tools for
    simulation, pre-processing, and image reconstruction, providing a multitude of
    different regularization options. In this work, we have used this toolbox for
    compressed sensing reconstructions of the velocities, using total variation in
    time as for the regularization. Examples of the reconstructed velocities, using
    different masks and acceleration factors, are shown in Figure \ref{fig:velocity_examples}.

\begin{figure}[!hbtp]
        \centering
        \subfloat[Spiral, $R=8$]{
        \includegraphics[trim=0 10 40 0, clip, width=0.3\textwidth]{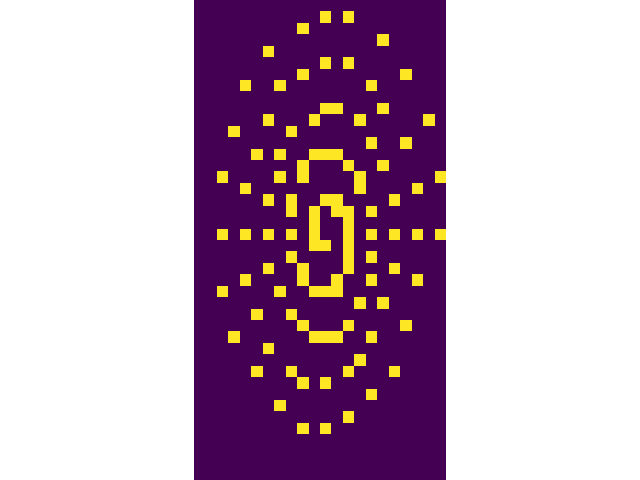} }
        \subfloat[Spiral, $R = 16$]{
        \includegraphics[trim=0 10 40 0, clip, width=0.3\textwidth]{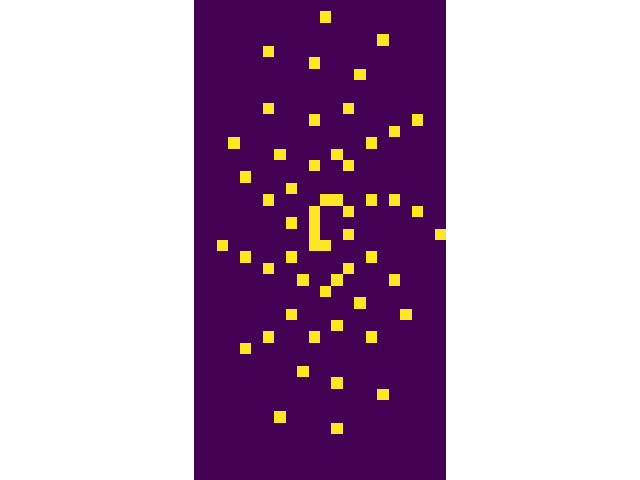} }
        \subfloat[Spiral, $R=32$]{
        \includegraphics[trim=0 10 40 0, clip, width=0.3\textwidth]{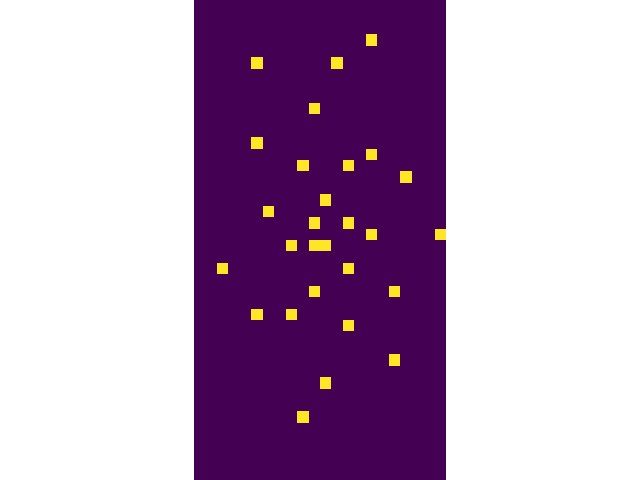} }\\
        \subfloat[Gaussian, $R=8$]{
        \includegraphics[trim=0 10 40 0, clip, width=0.3\textwidth]{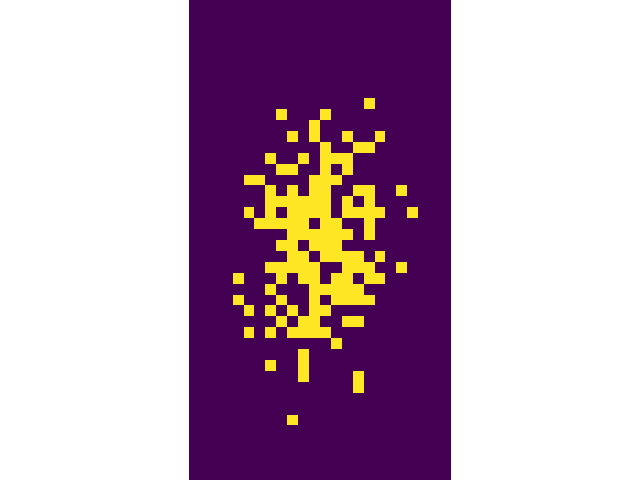} }
        \subfloat[Gaussian, $R=16$]{
        \includegraphics[trim=0 10 40 0, clip, width=0.3\textwidth]{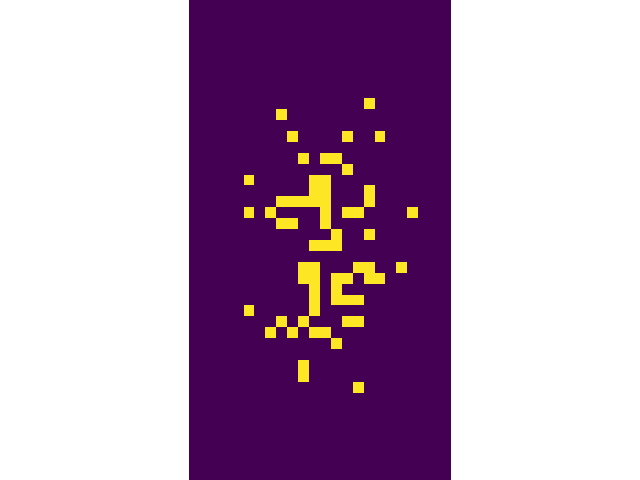} }
        \subfloat[Gaussian, $R=32$]{
        \includegraphics[trim=0 10 40 0, clip, width=0.3\textwidth]{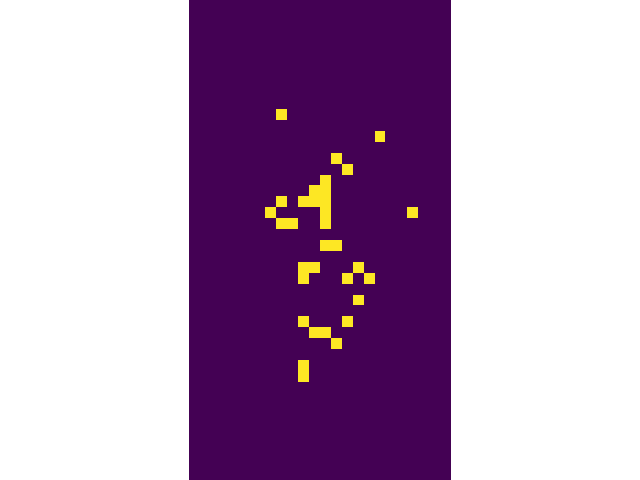} }\\

        \caption{Sampling masks}
        \label{fig:masks}
    \end{figure}

    Next, the sampling mask is applied to these simulated frequency space
    measurements. We take a 2D subsampled mask in the $x-y$-plane and sample fully in
    the $z$-direction as in \cite{peper2020highly}. We consider different subsampling
    rates $R = \frac{N_{sampled}}{N_{total}} = {8, 16, 32}$, with two different masks:
    the pseudo-spiral mask and the pseudo-random Gaussian mask, which is sampled
    according to a Gaussian probability distribution, as shown in Figure
    \ref{fig:masks}. For the pseudo-spiral mask, the points are placed evenly on a
    cartesian grid along a spiral with six turns and a final radius reaching the edge
    of the mask.

    Additionally, we require measurements of the magnitude itself. The value of the
    background phase is treated as known instead.

\begin{figure}[!hbtp]
        \centering
        \subfloat[K-space phase, fully sampled]{
        \includegraphics[trim=140 10 80 0, clip, width=0.3\textwidth]{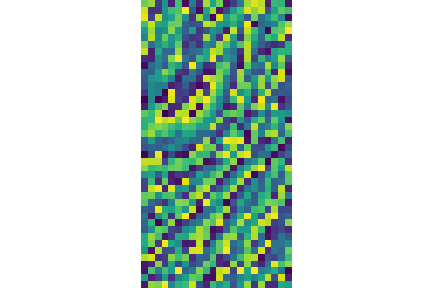} }
        \subfloat[True magnitude of the magnetization]{
        \includegraphics[trim=140 10 80 0, clip, width=0.3\textwidth]{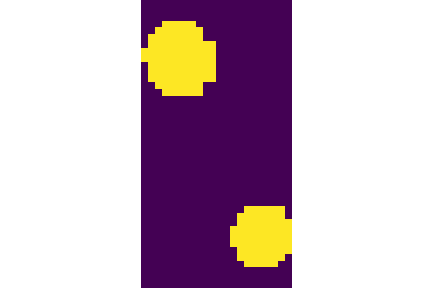} }
        \subfloat[True velocity in z-direction]{
        \includegraphics[trim=0 0 0 0, clip, width=0.3\textwidth]{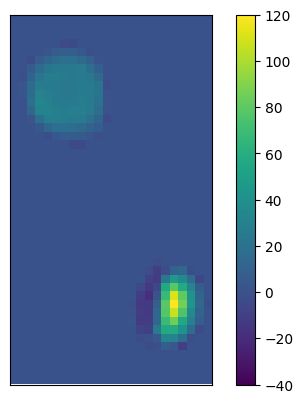} }\\
        \subfloat[K-space magnitude, fully sampled]{
        \includegraphics[trim=140 10 80 0, clip, width=0.3\textwidth]{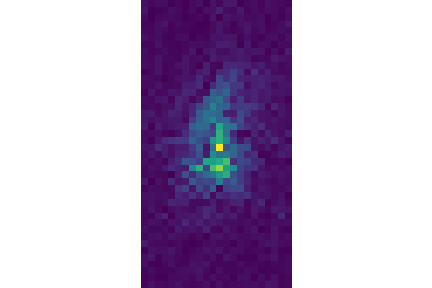} }
        \subfloat[Magnitude of the magnetization reconstructed with BART from $R=2$]{
        \includegraphics[trim=140 10 80 0, clip, width=0.3\textwidth]{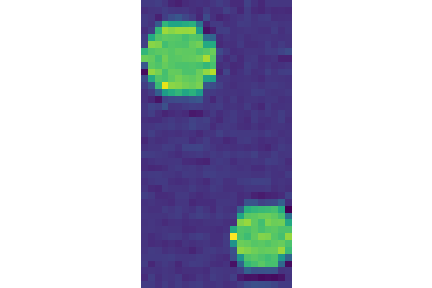} }
        \subfloat[Velocity in z-direction reconstructed with BART from R=2]{
        \includegraphics[trim=0 0 0 0, clip, width=0.3\textwidth]{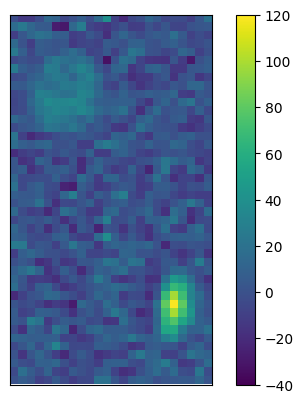} }\\

        \caption{Examples of simulated measurements, taken at a slice in the $z$-direction.}
        \label{fig:examples}
    \end{figure}

    As anatomical images are usually readily available and the magnitude generally
    does not depend on the encoding direction, making a cheap 2D/3D acquisition
    feasible, we consider a reconstructed magnitude from the measurements with $R = 2$
    and a Gaussian mask using temporal $l1$-regularization in BART with a
    regularization parameter $\lambda = 0.001$.
    Examples of the phase and magnitude of the measurements are shown in Figure
    \ref{fig:examples}.

    \subsubsection{Inverse problem setup}

    Without pressure measurements, not all the Windkessel resistances in the system
    can be uniquely determined at once. Therefore we have fixed the values for
    $\Gamma_4$ (the outlet in the descending aorta), and estimate  the amplitude of
    the inflow $U$ as well as the distal resistances $R_{d,k}$, $k=1,2,3$ of the
    remaining Windkessel boundary conditions.

    We consider two different initial guesses for these parameters, as in \cite{garay2022}:
    \begin{enumerate}
        \item "low" guess: $U = 40$, $R_{d,1}=4000$, $R_{d,2}=4000$, $R_{d,3}=4000$
        \item "high" guess: $U = 150$, $R_{d,1}=20000$, $R_{d,2}=20000$, $R_{d,3}=20000$
    \end{enumerate}
    while the target values are $U=75$, $R_{d,1}=7200$, $R_{d,2}=11520$, $R_{d,3}=11520$.

    The initial standard deviation for the reparameterized parameters was set to
    $0.5$, i.e.~$P^0 = 0.5\mathbb{I}$, meaning that the prior models that there is
    $\approx 95\%$ probability that the target value will lie within the range of
    half/twice the initial guess. The standard deviation of the noise was estimated
    from the initial time step by computing the standard deviation of
    \begin{equation*}
    \bs{Y}^0 - \bs{M}(t^0)\odot \exp\left(i\phi_{back}(t^0))\right)\odot \bs{S}
\end{equation*}
which is the noise of the data under the assumption that the velocity at the first
time step is zero. The estimated standard deviation for each case is shown in Table
\ref{tab:noise}. The estimation of the noise in this case is impacted by noise in the
reconstructed magnitude.

\begin{table}[!hbtp]
    \centering
    \begin{tabular}{|l|cc|}
        \hline
        Acceleration factor & Gaussian & Spiral\\
        \hline
        R1/actual value & \multicolumn{2}{c|}{15.526}\\
        R8 & 15.829 & 13.813 \\
        R16 & 18.716 & 13.791 \\
        R32 & 23.043 & 13.844\\
        \hline
    \end{tabular}
    \caption{Estimated standard deviation of the noise for each mask and acceleration
    factor, using the magnitude reconstructed from R2 for acceleration factors greater than 1}
    \label{tab:noise}
\end{table}

\subsection{Phantom data}

\subsubsection{Flow phantom measurements}
We used data from a flow experiment on a phantom of the carotid artery reported in
\cite{peper2020highly}. The phantom is made of a distensible silicon and is suspended
in water. It takes the shape of a bifurcating tube, to simulate a carotid artery, and
a backflow tube. A pump generates a pulsatile flow with a rate of roughy 60 bpm,
simulating a cardiac cycle lasting 1s.

\begin{figure}[!hbtp]
    \centering
    \subfloat[K-space magnitude of coil 4]{
    \includegraphics[trim=0 5 0 5, clip, width=0.25\textwidth]{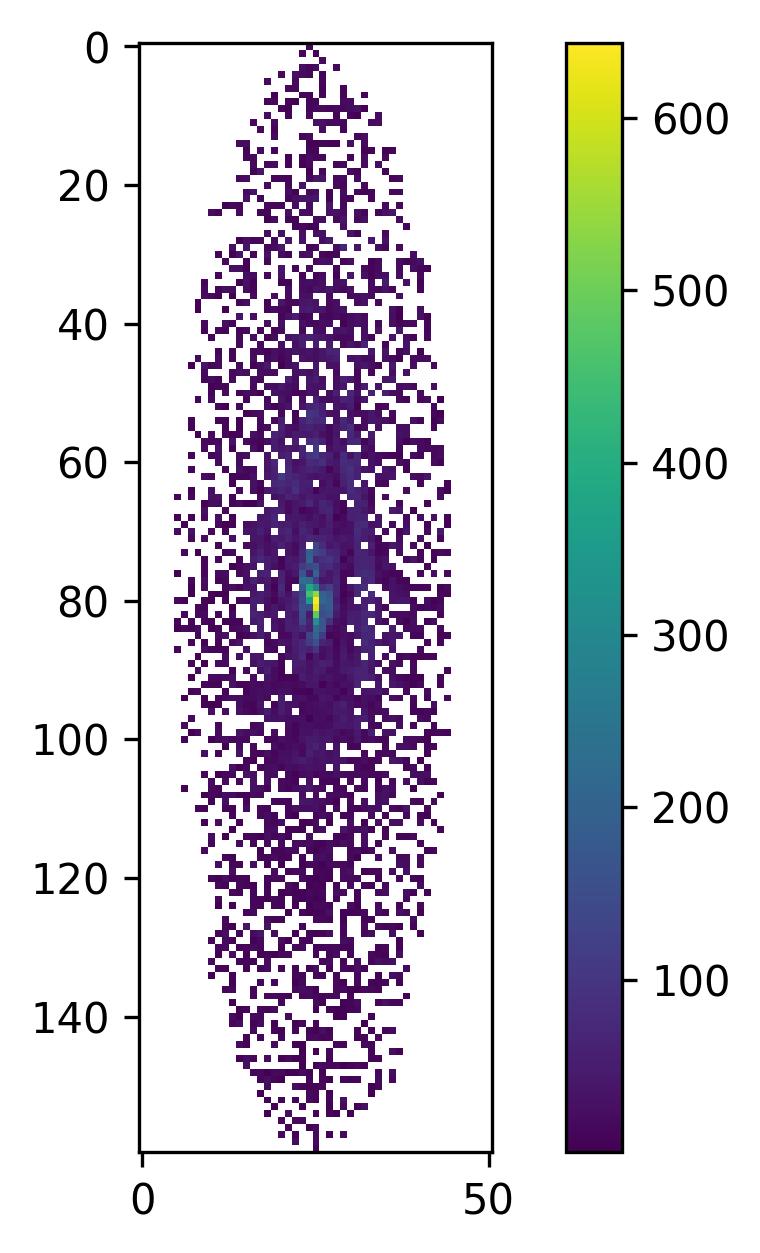} }
    \subfloat[x-component of the velocity reconstructed with BART from R=2]{
    \includegraphics[trim=0 5 0 5, clip, width=0.25\textwidth]{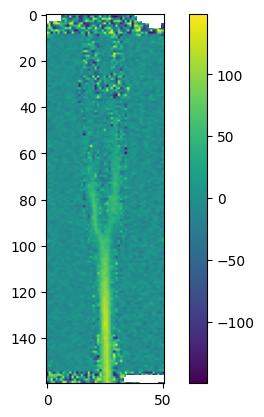} }
    \subfloat[Zero-filled magnitude of the magnetization of coil 4]{
    \includegraphics[trim=0 5 0 5, clip, width=0.25\textwidth]{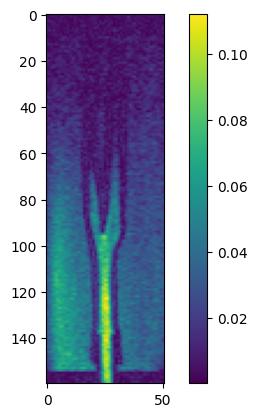} }
    \subfloat[Zero-filled magnitude of the magnetization of coil 6]{
    \includegraphics[trim=0 5 0 5, clip, width=0.25\textwidth]{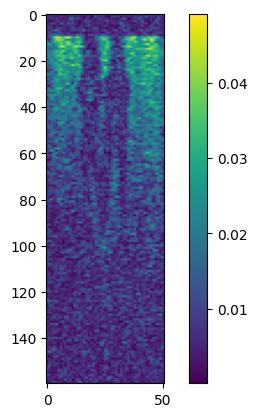} }\\

    \caption{Examples of the phantom measurements. Reconstructed magnitudes of the
    magnetization of different coils show the different coil sensitivities.}
    \label{fig:phantom_examples}
\end{figure}

The 4D Flow MRI scan was performed with a 3T Ingenia scanner by Philips Healthcare.
All three velocity directions were acquired, plus one acquisition with no encoding
gradient to acquire the background phase. The scan parameters were set to TR = 8.9ms,
TE = 4.5ms, FA = 8$^\circ$, $venc$=150cm/s. The matrix size was [160, 160, 40] with a
spatial resolution of [0.8mm, 0.8mm, 0.8mm] and a temporal resolution of 19 frames per
cardiac cycle. The scan was accelerated with an acceleration factor of R=2. The scan
used a 32-channel coil leading to 15 independent coil measurements, resulting in a
total matrix size of $[160, 160, 40, 4, 15, 19]$. Examples of the measurements, both
the frequency space magnitude and different reconstructions, are shown in Figure
\ref{fig:phantom_examples}.

\subsubsection{Forward problem setup}
\begin{figure}[!hbtp]
    \centering
    \subfloat[Segmented, structured mesh]{
    \includegraphics[trim=0 10 40 0, clip, width=0.8\textwidth]{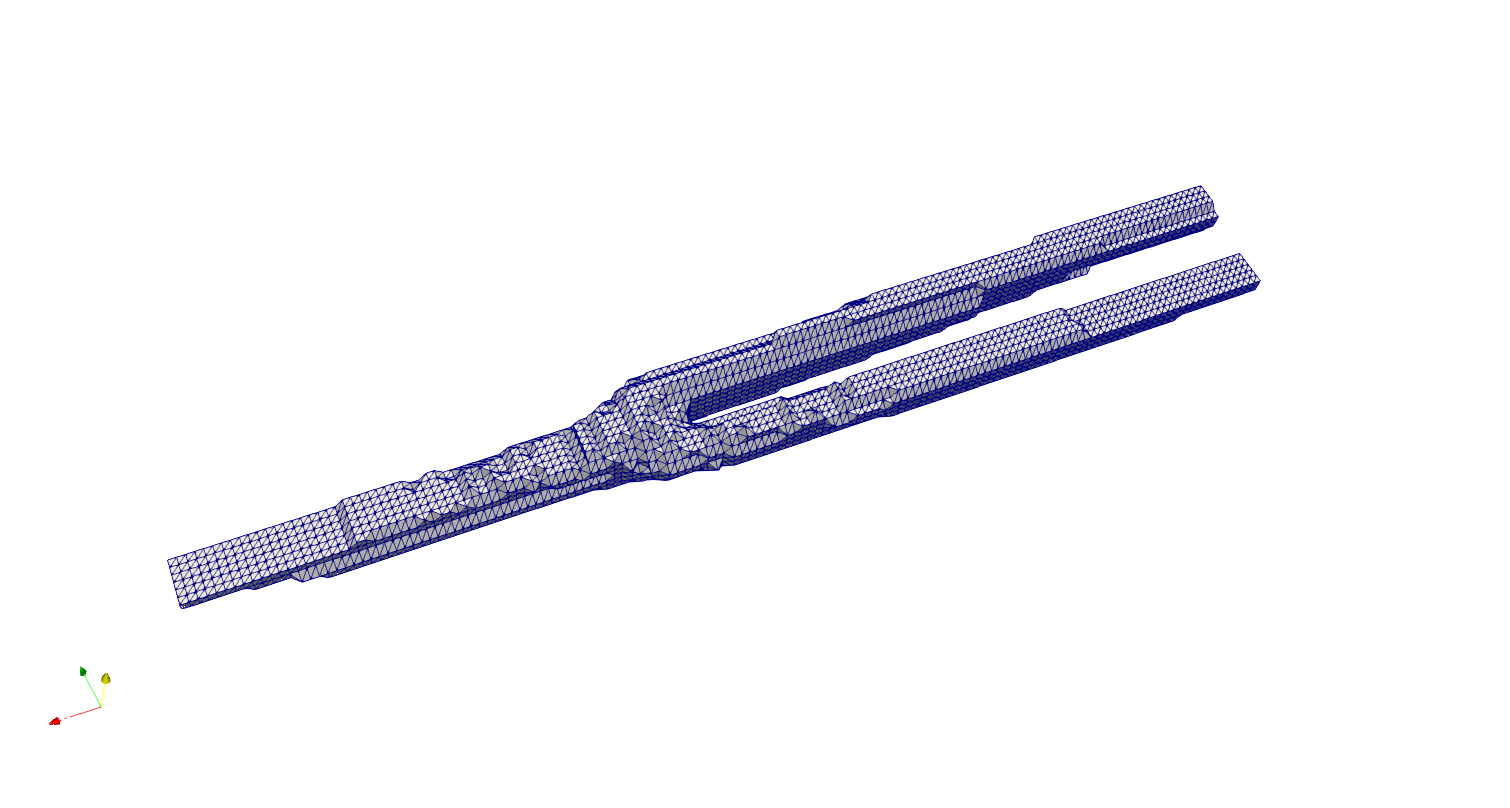} } \\
    \subfloat[Unstructured fine mesh]{
    \includegraphics[trim=0 10 40 0, clip, width=0.8\textwidth]{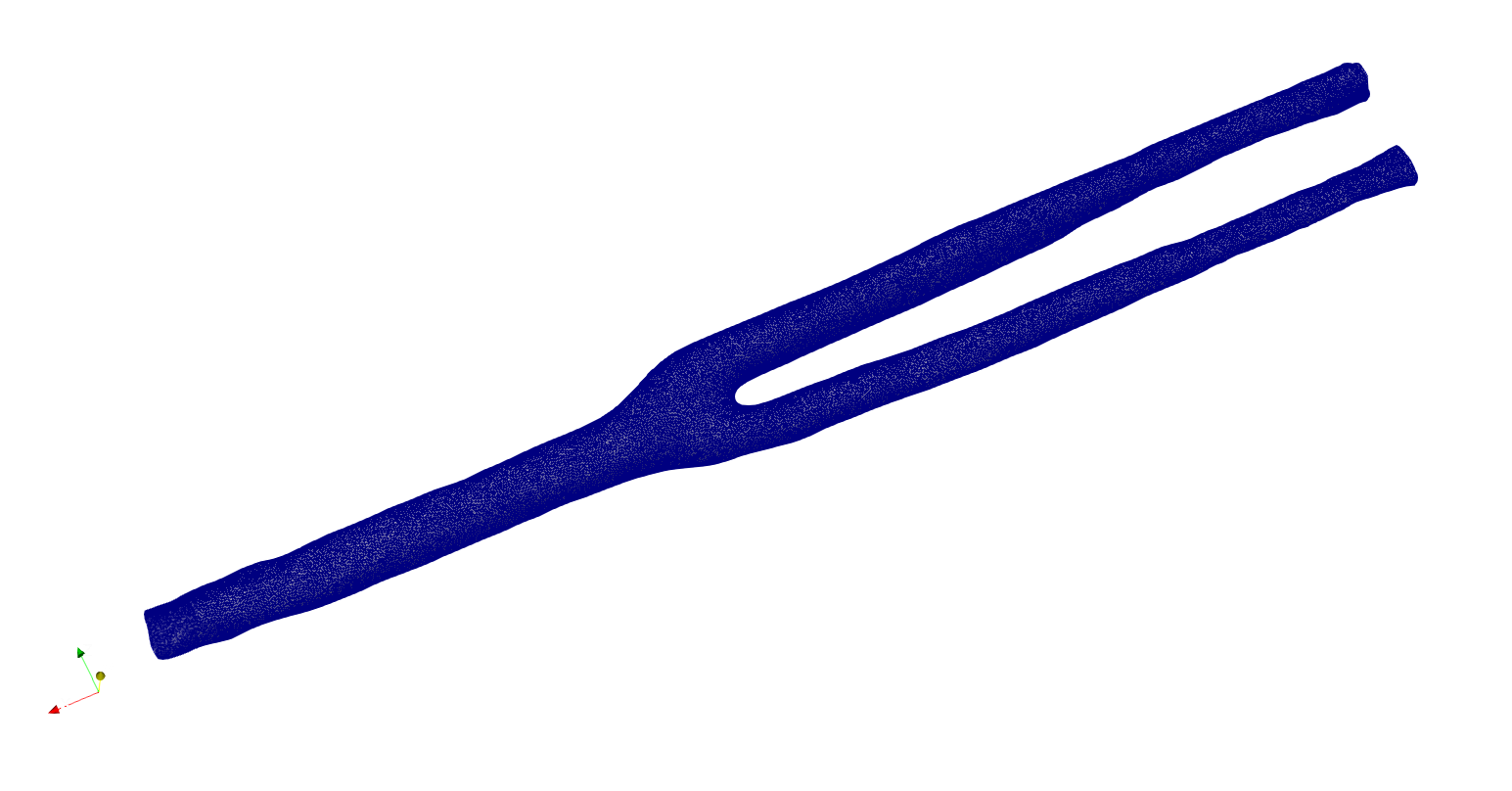} }

    \caption{Meshes of the carotid phantom}
    \label{fig:phantom_meshes}
\end{figure}

A structured mesh was created from a segmentation of the lumen, matching the spatial
resolution of the scan. From this an expanded, higher-resolution mesh was created
using Blender, Meshlab and gmsh. This mesh was used for the forward problem. The two
meshes are depicted in Figure \ref{fig:phantom_meshes}.

We again model the blood flow with the incompressible Navier-Stokes equations, with
the same choices for the physical parameters. The inflow boundary is modelled by a
Dirichlet boundary condition with an inflow
\begin{equation*}
    \bs{u}_{in} = \bs{u}_{profile}U f(t)
\end{equation*}
where $\bs{u}_{profile}$ is a flow profile provided by the solution of a Stokes
problem in the domain, $U$ is a constant amplitude, and $f(t)$ defines the temporal profile as
\begin{equation*}
    f(t) =
    \begin{cases}
        \sin(\frac{\pi t}{T}) \text{ if } t \leq \frac{3}{4} T \\
        sin(\frac{3}{4}\pi)(1-t+\frac{3}{4}T)\exp^{-(t-\frac{3}{4}T\beta}) \text{ if }
        t > \frac{3}{4} T
    \end{cases}
\end{equation*}
with $T = 0.64$, $\beta = 5$. The temporal profile was determined empirically to match
the shape of the flow rate in the phantom.

The two outlets were modelled with resistance boundary conditions
\begin{align*}
    P_l &= R_{p, l} Q_l \text{ on } \Gamma_l\\
    Q_l &= \int_{\Gamma_l} \bs{u}\cdot \bs{n} dx
\end{align*}
As the flow division between the outlets is determined by the proportion of the
resistances, we can fix $R_{p, 1} = 100$ at an arbitrary value.

\subsubsection{Inverse problem setup}
We estimate the inflow amplitude $U$ and the resistance $R_{p, 2}$ of the boundary
condition at the right outlet. The initial guess is $U = 100$, $R_{p,2} = 100$,
i.e.~an equal division of flow between the outlets. The initial standard deviation was
set to $0.5$.

We are reparameterizing such that $\bs{\theta} = \bs{\theta}^0  \odot \bs{\nu}$ with
$\bs{\theta}^0$ the initial guess for the parameters and the filter being applied to
$\bs{\nu}$, with the initial value of $\bs{\nu} = \bs{1}$. This was done because the
exponential reparameterization used before proved unstable in this case.

To extend the filter to multiple coils, each of the coil measurements was treated as a
separate, independent measurement, resulting in 15 measurements per time step. The
variance of the noise for each coil was estimated from the first measurement by
computing the standard deviation of
\begin{equation*}
\bs{Y}^0 - \bs{M}(t^0)\odot \exp\left(i\bs{\phi}_{back}(t^0))\right)\odot \bs{S}
\end{equation*}
as before, which neglects potential spatial variation of the variance due to the
sensitivity of the coils.

The observation operator requires measurements of the magnitude of the magnetization
and the background phase. The background phase is reconstructed as
\begin{equation}
\bs{\phi_{back}} = \angle \mathcal{F}^{-1}(\bs{Y}_{back})
\end{equation}
where $\bs{Y}_{back}$ is a zero-filled measurement acquired with $R=2$ and no encoding
gradient, i.e.~capturing only the background magnetization and no fluid velocity.

For the magnitude, we consider two options. The first option is
\begin{equation}
\bs{M} = |\mathcal{F}^{-1}(\bs{Y})|
\end{equation}
for each velocity direction, where $\bs{Y}$ is a zero-filled measurement acquired with
$R=2$ in that velocity direction. This however assumes that a highly sampled
measurement is available. Therefore the other option is to use
\begin{equation}
\bs{M} = |\mathcal{F}^{-1}(\bs{Y}_{back})|
\end{equation}
for all velocity directions, using the same measurement as for $\phi_{back}$. This
would only require a highly sampled measurement of one out of four encoding gradients.

\begin{figure}[!hbtp]
\centering
\subfloat[Original, $R = 2$]{
\includegraphics[trim=120 10 120 30, clip, width=0.3\textwidth]{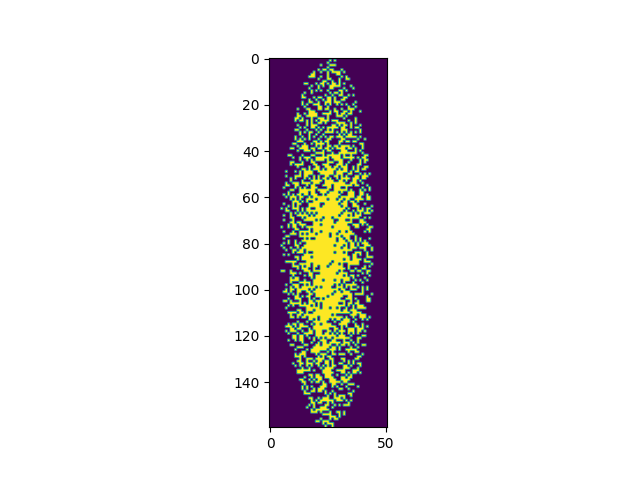} } \\

\subfloat[Spiral, $R = 16$]{
\includegraphics[trim=140 10 140 30, clip, width=0.23\textwidth]{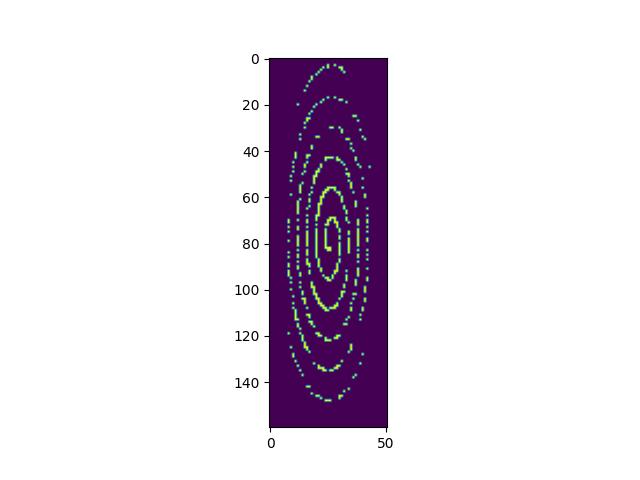} }
\subfloat[Spiral, $R=32$]{
\includegraphics[trim=140 10 140 30, clip, width=0.23\textwidth]{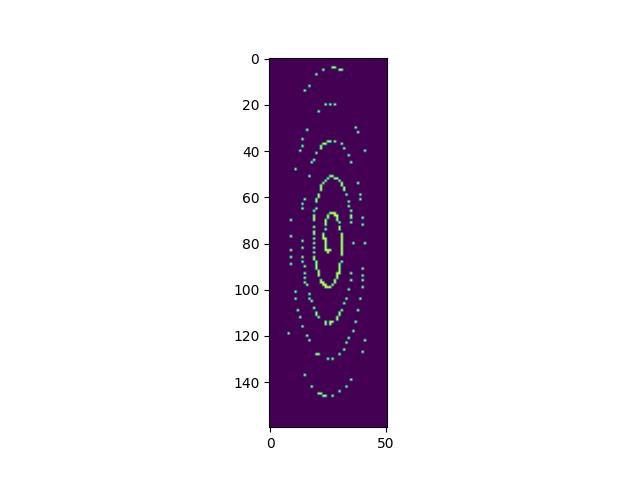} }
\subfloat[Spiral, $R=64$]{
\includegraphics[trim=140 10 140 30, clip, width=0.23\textwidth]{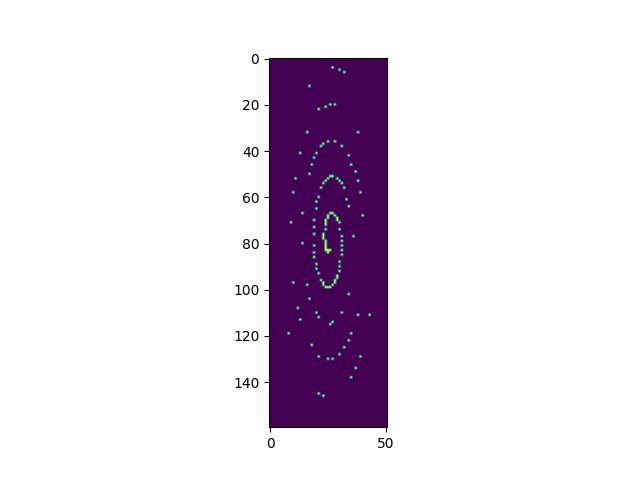} }
\subfloat[Spiral, $R=128$]{
\includegraphics[trim=140 10 140 30, clip, width=0.23\textwidth]{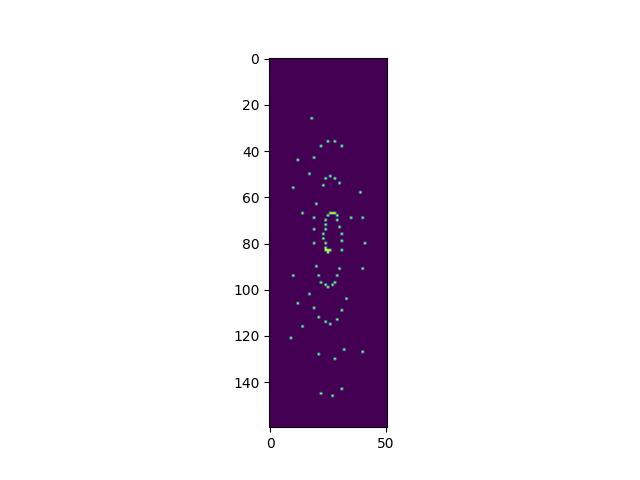} }\\
\subfloat[Gaussian, $R=16$]{
\includegraphics[trim=140 10 140 30, clip, width=0.23\textwidth]{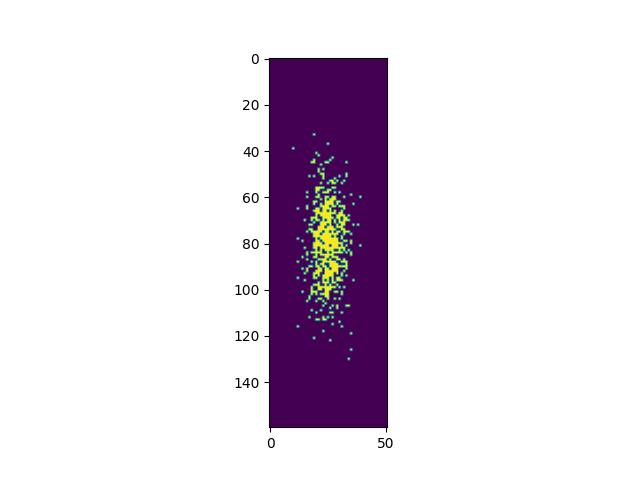} }
\subfloat[Gaussian, $R=32$]{
\includegraphics[trim=140 10 140 30, clip, width=0.23\textwidth]{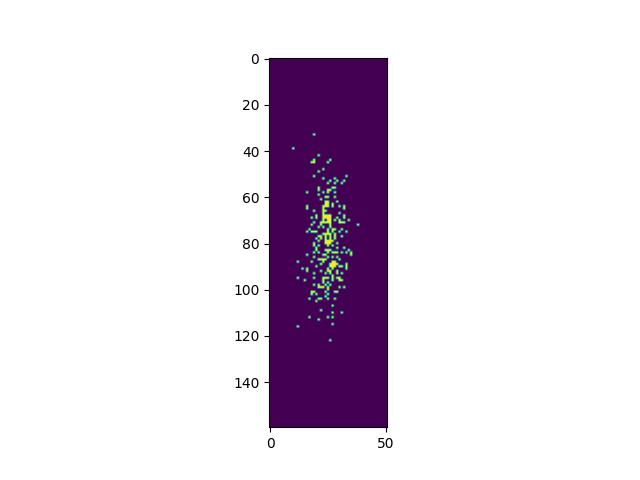} }
\subfloat[Gaussian, $R=64$]{
\includegraphics[trim=140 10 140 30, clip, width=0.23\textwidth]{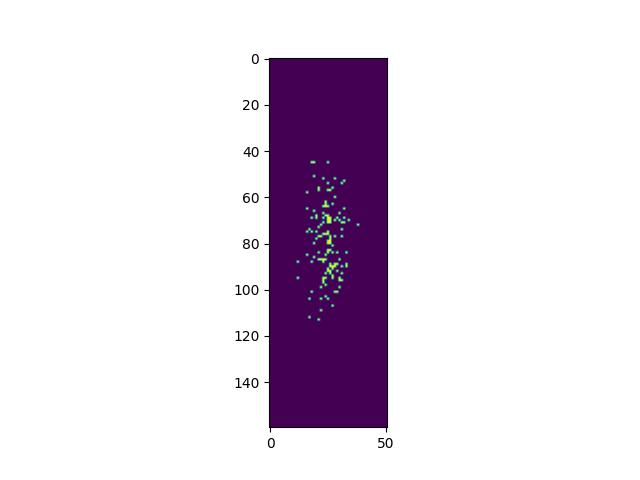} }
\subfloat[Gaussian, $R=128$]{
\includegraphics[trim=140 10 140 30, clip, width=0.23\textwidth]{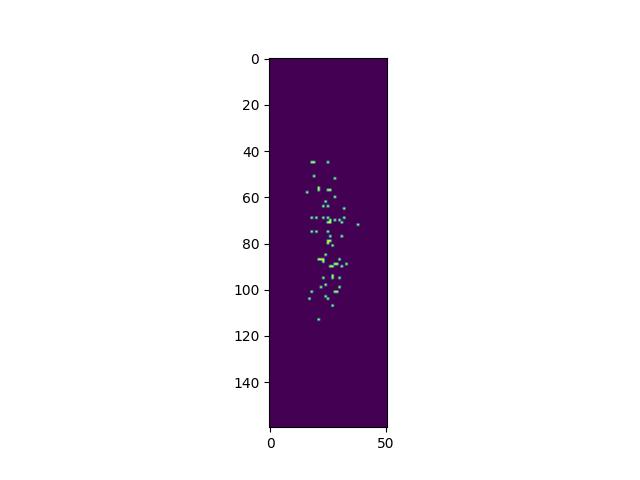} }
\caption{Sampling masks for the phantom data, taken at the $k_x = 0$ slice at time step 0}
\label{fig:phantom_masks}
\end{figure}

The original data with an acceleration factor of 2 were acquired with an incoherent
pseudo-spiral sampling in the $k_y-k_z$ direction with a fully sampled
$k_x$-direction. We further undersample this by applying masks with a Gaussian or
spiral sampling pattern onto the sampling mask for $R=2$ to achieve higher
acceleration factors $R = 16, 32, 64, 128$. In contrast to the masks for the synthetic
data, this results in different masks for each time step and velocity direction.
Examples of the resulting masks are shown in Figure \ref{fig:phantom_masks}.

The measurements are assumed to be placed at every $0.053s (\approx 1s/19)$. For the
inverse problem, the first two measurements (at $0.053s$ and $0.106s$) are omitted as
the presence of negative velocity values in these measurements could interfere with
the Kalman filter, since this could lead to negative values in the particles which
would be unphyscial for the forward simulation.

\section{Results}

\label{sec:results}

\subsection{Synthetic data}
\label{sec:synth_results}

\noindent\textit{Computational effort.} The parameter estimation with either velocity
data or frequency data requires the same amount of CPU time (around 13 minutes using
an AMD 7763 CPU), though it can be run in parallel to reduce the walltime. Using
velocity measurements uses slightly less memory however, as no complex-value data have
to be saved (888.6MB compared to 959MB for the frequency-space data). Performing the
reconstruction of the velocity data with BART requires an additional 1:30min of CPU
time and 660MB of memory.

\bigskip
We compare the results of our method (directly estimating from frequency-space
measurements) to estimating parameters from velocity measurements reconstructed with
compressed sensing using BART with a $l1$-regularization in time. The regularization
parameter was determined empirically by visual evaluation for each combination of
acceleration factor and mask and is listed in Table \ref{tab:lambda}. For the
estimation, we are using the method presented in \cite{garay2022} to take advantage of
knowing the $R=2$ magnitude.

\bigskip
\noindent\textit{Comparing the reconstruted flows.} In order to accurately reflect the
different impact of the parameters on the flow, we compute the error in terms of the
flow computed from the estimated parameters. Using the estimated parameters, we solve
the forward problem \ref{eq:forward} again with a time step of $dt = 0.001$. The error
is calculated as
\begin{equation}
e = \frac{||\bs{u}_{ref} - \bs{u}_{recon}||_{2}}{||\bs{u}_{ref}||_2}
\end{equation}\label{eq:error}
where $\bs{u}_{recon}$ and $\bs{u}_{ref}$ are vectors consisting of the velocity
values for all components at each point in the geometry at each point in time stacked
together, using the estimated parameter values for $\bs{u}_{recon}$ and the true
parameter values for $\bs{u}_{ref}$.

\begin{table}[!hbtp]
\centering
\begin{tabular}{|l|r|r|r|}
    \hline
    Acc. Factor & Spiral & Gaussian\\
    \hline
    $R1$ & \multicolumn{2}{c|}{$0.0001$} \\
    $R8$   & $0.1$ & $0.01$\\
    $R16$  & $1.0$ & $0.01$\\
    $R32$  & $10.0$ & $1.0$\\
    \hline
\end{tabular}
\caption{Regularization values $\lambda$ used for the reconstruction with BART.}
\label{tab:lambda}
\end{table}

\begin{figure}[!hbtp]
\centering
\subfloat[low initial guess]{
\includegraphics[trim=20 10 40 0, clip, width=0.48\textwidth]{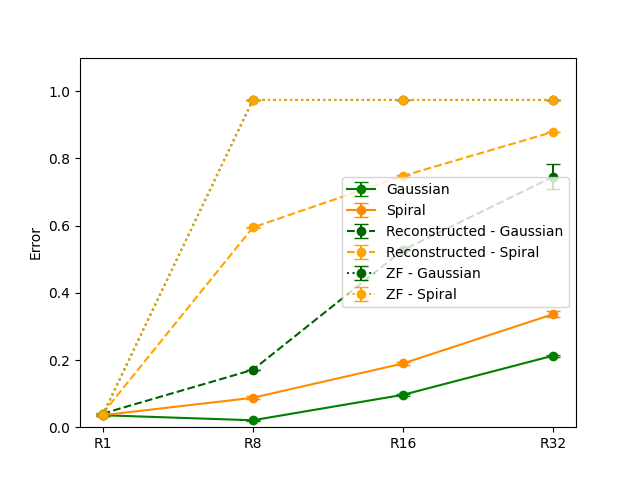} }
\subfloat[high initial guess]{
\includegraphics[trim=20 10 40 0, clip, width=0.48\textwidth]{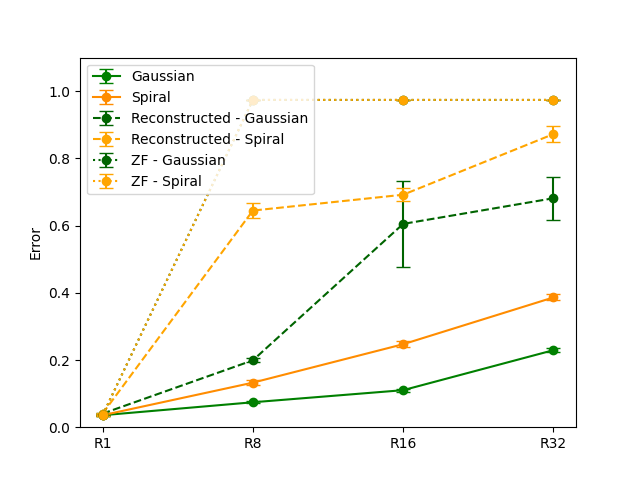} }\\

\caption{Error values for different acceleration factors and different masks. Dotted
    lines are the error values of the flow reconstructed from the parameters estimated by
    the inverse problem using zero-filled measurements, dashed lines for the flow
    reconstructed from the inverse problem from BART velocities, solid lines are from the
    inverse problem in frequency space. Low initial guess on the left, high initial guess
on the right. The bars indicate the standard deviation of the error.}
\label{fig:error_lines}
\end{figure}

The error values for all different acceleration factors for the low and high initial
guess are depicted in Figure \ref{fig:error_lines}, for our method as well as for
estimating the parameters using velocities reconstructed with BART, as well as the
error values of the velocity reconstructed with BART itself (without first estimating
the parameters). In the latter case, $u_{ref}$ was interpolated into the coarser mesh
to match the resolution of the reconstructed velocities. In all cases, estimating the
parameters and then reconstructing the flow achieves better results than using
zero-filled velocity measurements. For the k-space cost function, the errors of both
masks are very close to the error with fully sampled data for $R=8$, and increase for
$R=16$ and $R=32$. The Gaussian mask achieves lower error values for all three
acceleration factors.

Considering the velocity measurements from data reconstructed using BART, it can be
seen that the Gaussian mask performs considerably better than the spiral mask,
especially for $R=8$. This matches the expectation, as a pseudo-random masks leads to
incoherent artifacts, which can be easily excluded with a temporal regularizer.
Nonetheless, the error increases drastically for $R=16$ and $R=32$. The spiral mask
shows a strong increase in error for $R=8$ already and remains high with higher
subsampling rates.

Using frequency measurements directly outperforms using the reconstructed velocity
measurements for all acceleration factors except for the fully sampled case.

Both the high and the low initial guess show the same pattern. As such, from here on,
we will limit ourselves to showing results for the low initial guess as the ones for
the high initial guess provide no additional information.

\begin{figure}[!hbtp]
\centering
\subfloat[$U$]{
\includegraphics[trim=20 10 40 0, clip, width=0.48\textwidth]{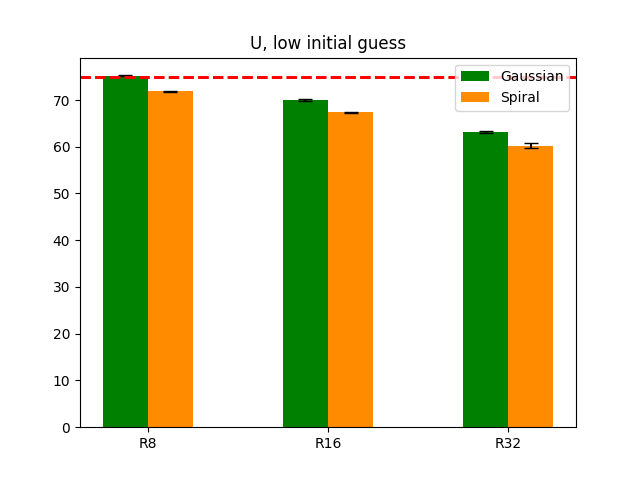} }
\subfloat[$R_{d, 1}$]{
\includegraphics[trim=20 10 40 0, clip, width=0.48\textwidth]{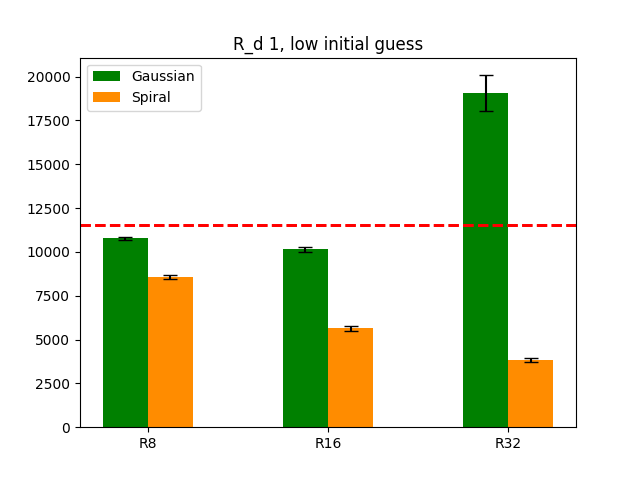} }\\
\subfloat[$R_{d, 2}$]{
\includegraphics[trim=20 10 40 0, clip, width=0.48\textwidth]{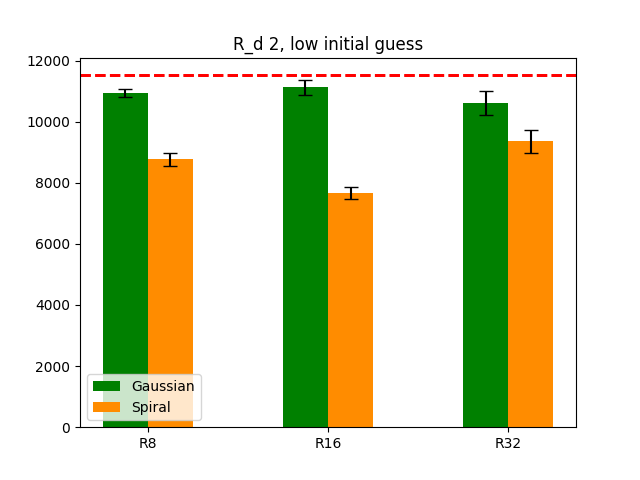} }
\subfloat[$R_{d, 3}$]{
\includegraphics[trim=20 10 40 0, clip, width=0.48\textwidth]{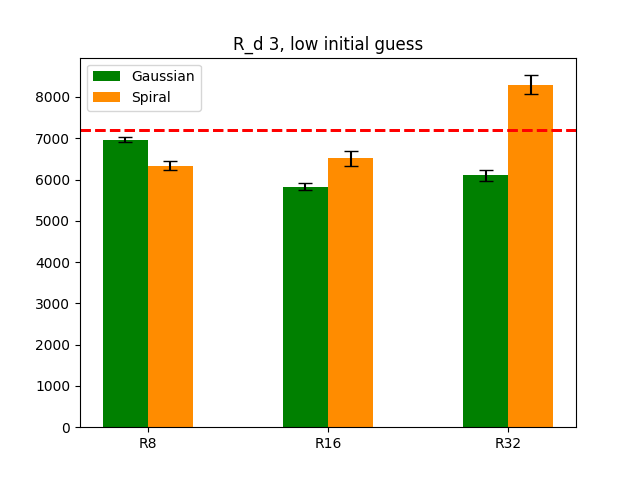} }\\

\caption{Estimated values of individual parameters. The dashed line indicates the true value.}
\label{fig:error_by_param}
\end{figure}

\bigskip
\noindent\textit{Comparing the parameter values.} By comparing the estimated parameter
values of each parameter individually, as seen in Figure \ref{fig:error_by_param}, it
is apparent that the success of the estimation differs depending on the parameter.
The inflow $U$ is estimated relatively accurately by all masks, with only a small
decrease with increasing $R$. In comparison, the distal resistance of the first
Windkessel outlet ($R_{d,1}$) is underestimated considerably and increasingly by the
spiral mask, while the gaussian mask estimates it well for $R=8$ and $R=16$, but
overestimates it severely for $R = 32$. $R_{d, 2}$ and $R_{d, 3}$ show the same
pattern, with the spiral mask underestimating the values compared to the Gaussian
mask, though more severely for $R_{d,2}$ than $R_{d, 3}$. For $R=32$, the spiral mask
slightly overestimates the value of $R_{d, 3}$.

\begin{figure}[!hbtp]
\centering
{
\includegraphics[trim=0 10 40 0, clip, width=0.7\textwidth]{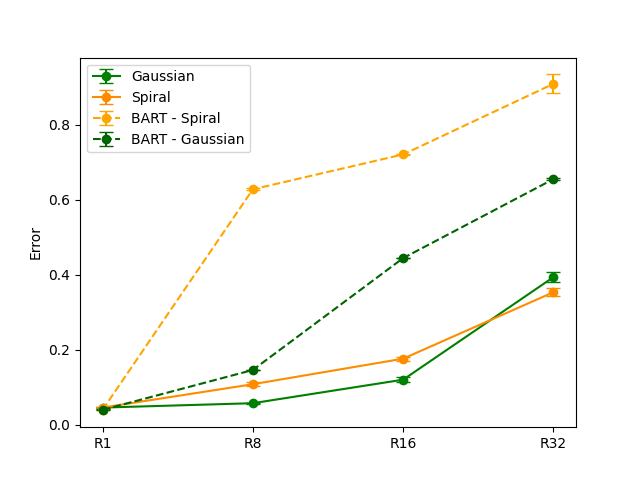} }

\caption{Error values for different acceleration factors using only the $z$-component
of the velocity.}
\label{fig:z}
\end{figure}

\bigskip
\noindent\textit{Using only a single velocity component.} We also consider the case
where measurements are only available for one of the velocity components, in this case
the $z$-component, which is equivalent to the foot-head direction in this setup.
Figure \ref{fig:z} shows the error values for this case. The error values show the
same patterns as when using all velocity components, and are only slightly higher.
Here there is a lesser difference in the error values between the spiral and Gaussian
mask for $R= 8$ and $R = 16$, and the spiral mask shows a flatter error curve, whereas
the error for the Gaussian mask increases with $R=32$. Again, using the frequency
measurements achieves lower error values than BART measurements for all subsampled data.

\begin{figure}[!hbtp]
\centering
\subfloat[$U$]{
\includegraphics[trim=0 10 40 0, clip, width=0.48\textwidth]{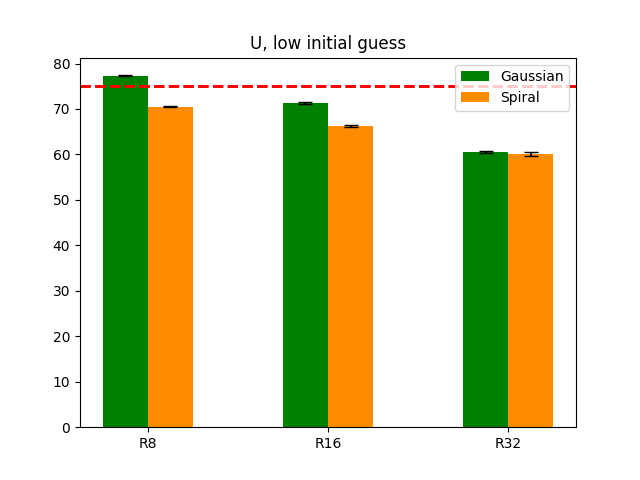} }
\subfloat[$R_{d, 1}$]{
\includegraphics[trim=0 10 40 0, clip, width=0.48\textwidth]{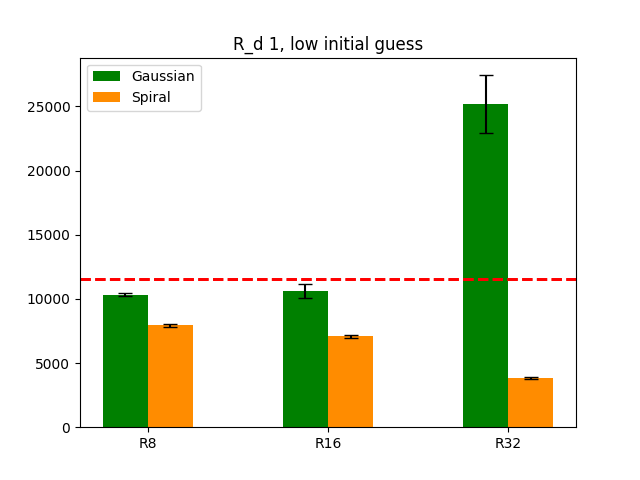} }\\
\subfloat[$R_{d, 2}$]{
\includegraphics[trim=0 10 40 0, clip, width=0.48\textwidth]{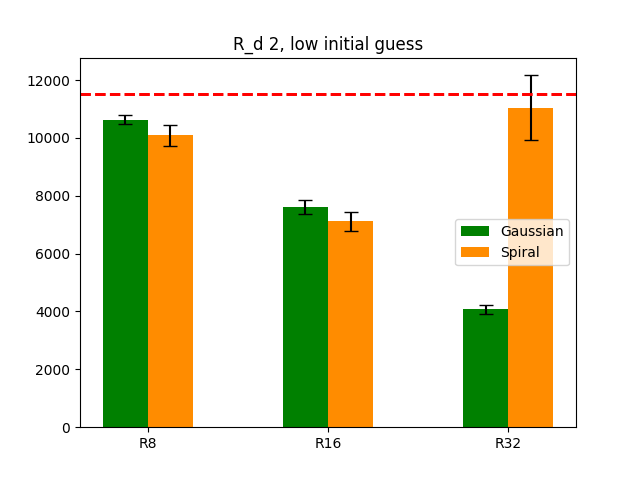} }
\subfloat[$R_{d, 3}$]{
\includegraphics[trim=0 10 40 0, clip, width=0.48\textwidth]{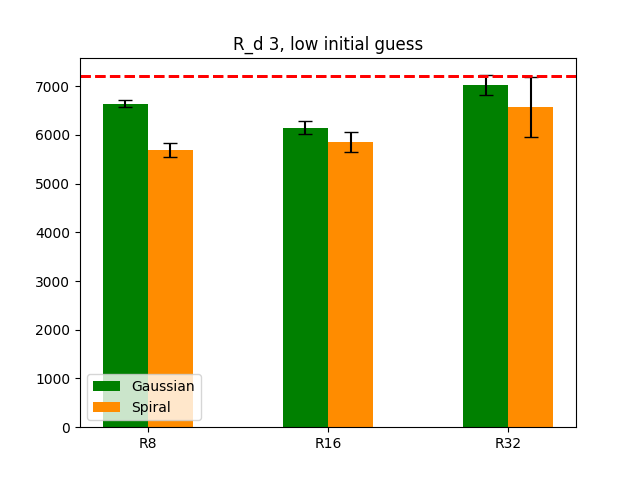} }\\

\caption{Estimated values of individual parameters using only the $z$-component of the
velocity. The dashed line indicates the true value.}
\label{fig:error_by_param_z}
\end{figure}
When considering each parameter separately in Figure \ref{fig:error_by_param_z}, the
same patterns and similar values persist as for all velocity components. The exception
is the Gaussian mask for $R_{d,2}$, which now significantly underestimates the values,
unlike the spiral mask.

\subsection*{Robustness to the choice of $venc$}

As described previously in Remark 2, a low $venc$ leads to a high signal-to-noise
ratio, but can lead to aliasing artifacts if the actual maximal velocity exceeds the
$venc$. The cost function used in \cite{garay2022} remedies this by distinguishing the
actual velocity from the wrapped ones when the physical parameters affect the velocity
on several voxels simultaneously. In the present work, the cost function corresponds
to a similar cost function as in  \cite{garay2022} -- but with an additional Fourier
transform -- therefore also including  the aliasing compensation.

\begin{figure}[!hbtp]
\centering
\subfloat[Spiral]{
    {
\includegraphics[trim=10 10 40 10, clip, width=0.47\textwidth]{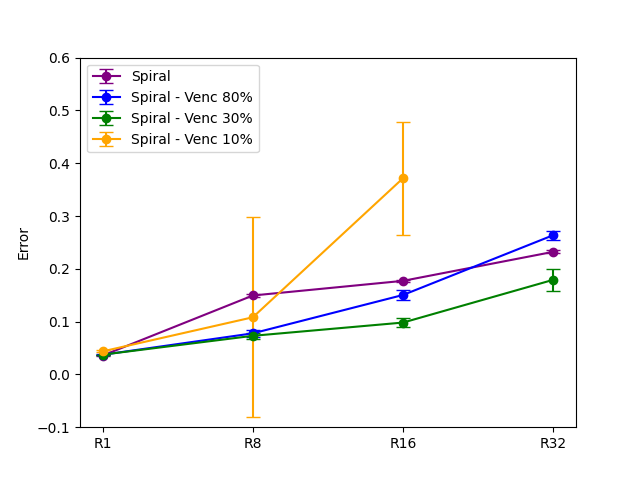} } }
\subfloat[Gaussian]{
    \includegraphics[trim=10 10 40 10, clip, width=0.47\textwidth]{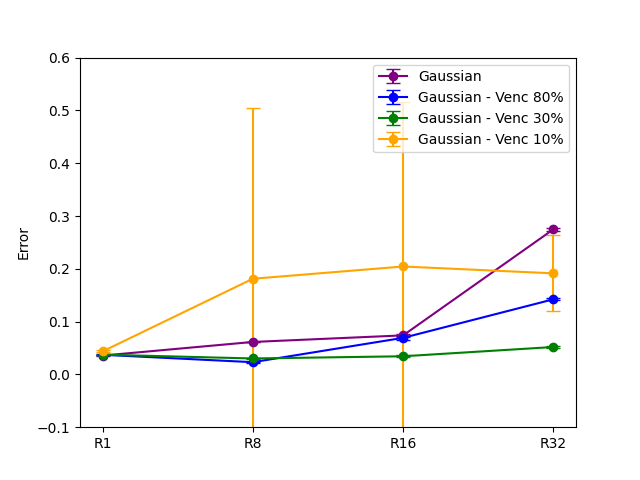}
}\\
\subfloat[BART with spiral mask]{
    {
\includegraphics[trim=10 10 40 10, clip, width=0.47\textwidth]{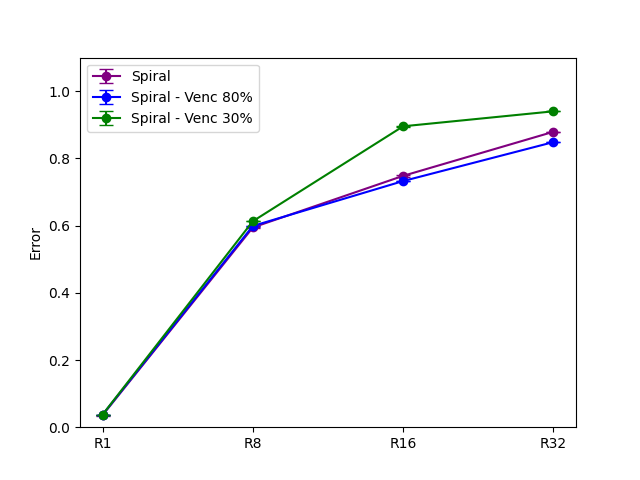} } }
\subfloat[BART with gaussian mask]{
\includegraphics[trim=10 10 40 10, clip, width=0.47\textwidth]{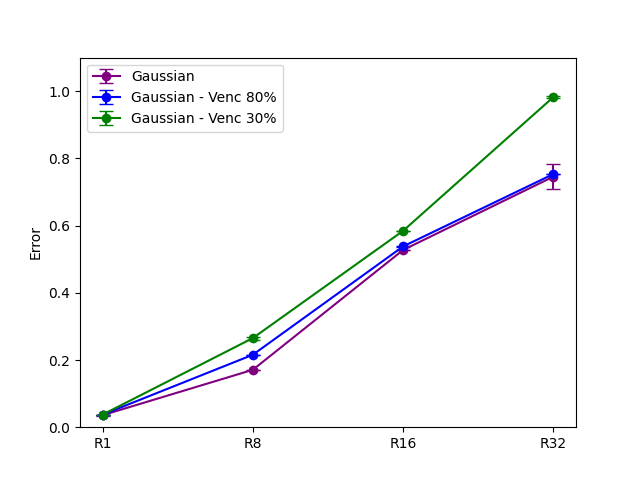} }

\caption{Error values for different values of the $venc$. Spiral mask on the left,
    Gaussian mask on the right. The value for $R32$ for the spiral mask for the lowest
    venc is excluded as the inverse problem failed to converge. For the BART
reconstructions, $venc = 10\%$ was not considered.}
\label{fig:venc_comparison}
\end{figure}

To investigate the robustness of our method to choosing $venc$s lower than the maximal
velocity, thus utilizing the higher sensitivity to the velocity in the signal, we
compare three additional $venc$ values corresponding to $80\%$, $30\%$, and $10\%$ of
the maximum velocity. The results are shown in Figure \ref{fig:venc_comparison}.

The results improve considerably with a decreasing $venc$ for 80\% and 30\%,
especially for higher acceleration factors. However, the smallest $venc$ of 10\% leads
to increased error values and standard deviations, with the Kalman filter even failing
to provide results for the spiral mask for $R=32$. The results for both masks perform
similarly with the decreasing $venc$, indicating that this response is inherent to the
inverse problem formulation -- as shown in \cite{garay2022} -- rather than a result of
the choice of mask.

In comparison, with the BART reconstructions, the error remains similar for $venc =
80\% V_{max}$ and then increases significantly for $venc = 30\% V_{max}$. This
behaviour can be explained from the fact that the reconstructed velocities were
already unwrapped prior to the parameter estimation as described in \cite{garay2022}.

\subsection*{Using an estimated magnitude of the magnetization}
In the prior results, we have been using a magnitude reconstructed from k-space
subsampled with $R=2$. To check the accuracy of this substitution, we now compare to
providing perfect knowledge of the magnitude. We also compare to using a constant
magnitude with a value of $0.5$, to model having little-to-no knowledge of the magnitude.

\begin{figure}[!hbtp]
\centering
\subfloat[Spiral]{
    {
\includegraphics[trim=10 10 40 10, clip, width=0.47\textwidth]{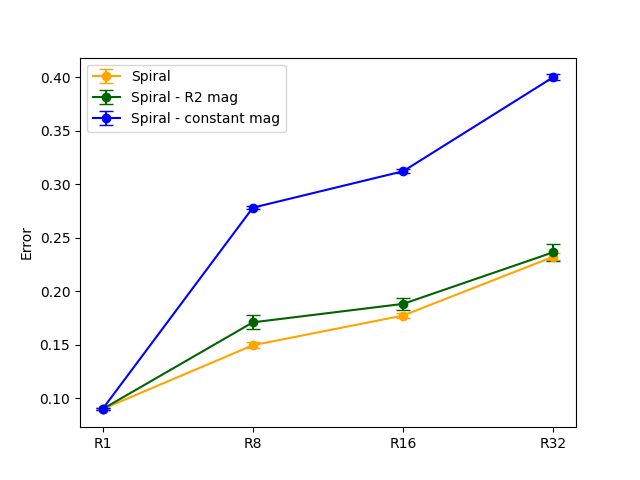} } }
\subfloat[Gaussian]{
    \includegraphics[trim=10 10 40 10, clip, width=0.47\textwidth]{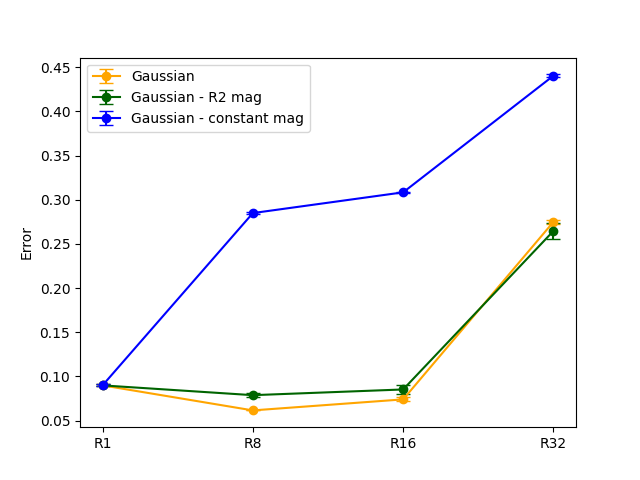}
}

\caption{Error values for different estimations of the magnitude}
\label{fig:mags}
\end{figure}

As can be seen in Figure \ref{fig:mags}, using the magnitude reconstructed from
frequency data with an acceleration factor of $R=2$ achieves results that are very
close to using the perfect magnitude. On the other hand, using a constant magnitude
leads to a significant increase in error that appears constant across the different
acceleration factors. In this case, there are also less differences between the two
masks. Nonetheless, for higher acceleration factors, even with a constant magnitude
the error is less than with velocity measurements reconstructed from BART.

\begin{figure}[!hbtp]
\centering
\subfloat[Spiral, $U$]{
\includegraphics[trim=0 10 0 0, clip, width=0.3\textwidth]{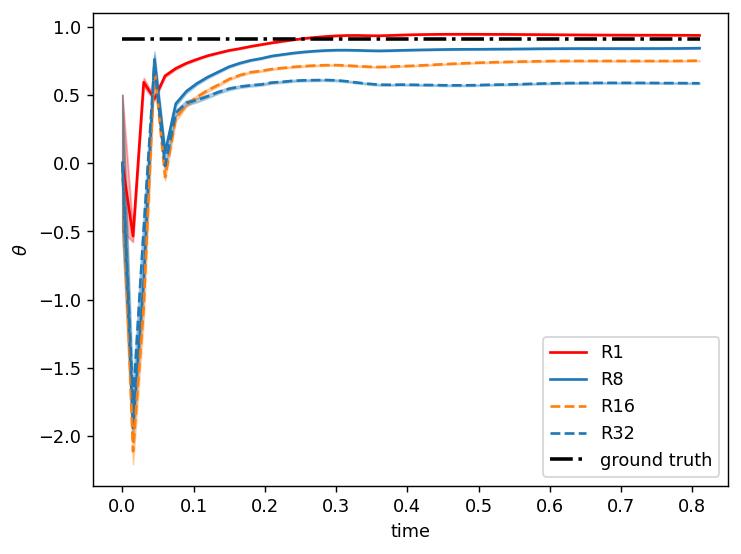} }
\subfloat[Gaussian, $U$]{
\includegraphics[trim=0 10 0 0, clip, width=0.3\textwidth]{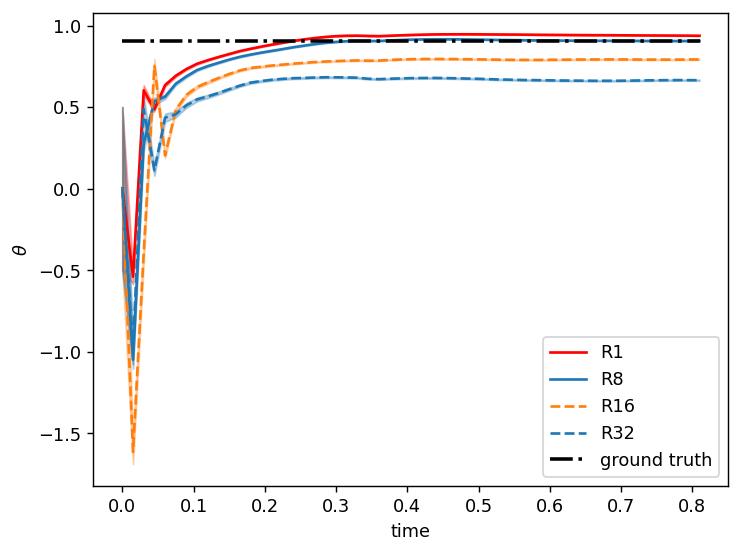} }
\subfloat[BART (Gaussian mask), $U$]{
\includegraphics[trim=0 10 0 0, clip, width=0.3\textwidth]{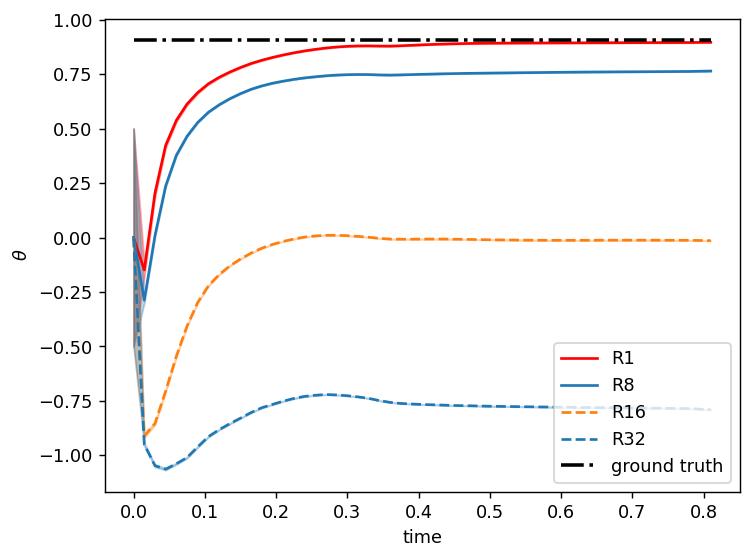} }\\
\subfloat[Spiral, $R_{d, 1}$]{
\includegraphics[trim=0 10 0 0, clip, width=0.3\textwidth]{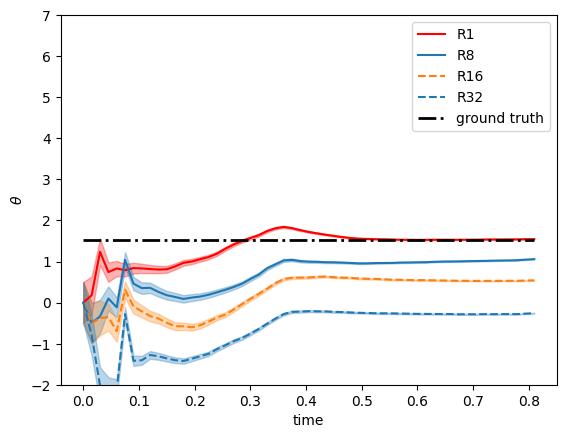} }
\subfloat[Gaussian, $R_{d, 1}$]{
\includegraphics[trim=0 10 0 0, clip, width=0.3\textwidth]{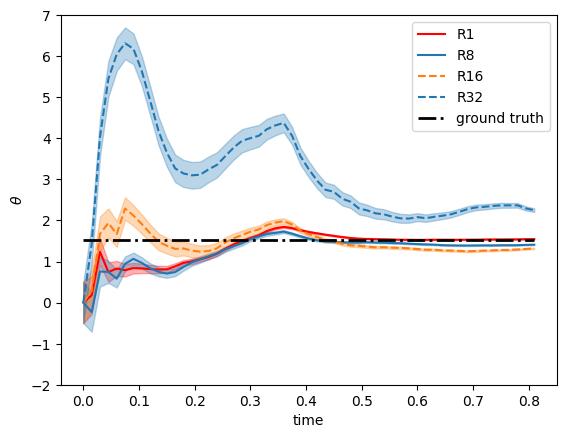} }
\subfloat[BART (Gaussian mask, $R_{d, 1}$)]{
\includegraphics[trim=0 10 0 0, clip, width=0.3\textwidth]{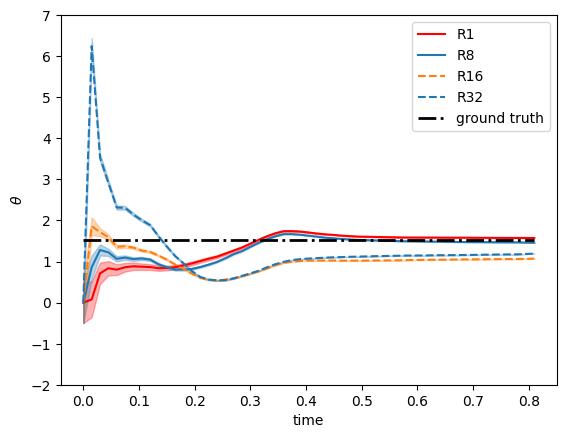} }\\
\subfloat[Spiral, $R_{d, 2}$]{
\includegraphics[trim=0 10 0 0, clip, width=0.3\textwidth]{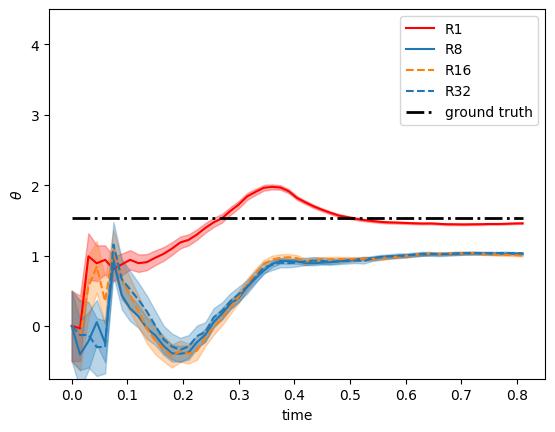} }
\subfloat[Gaussian, $R_{d, 2}$]{
\includegraphics[trim=0 10 0 0, clip, width=0.3\textwidth]{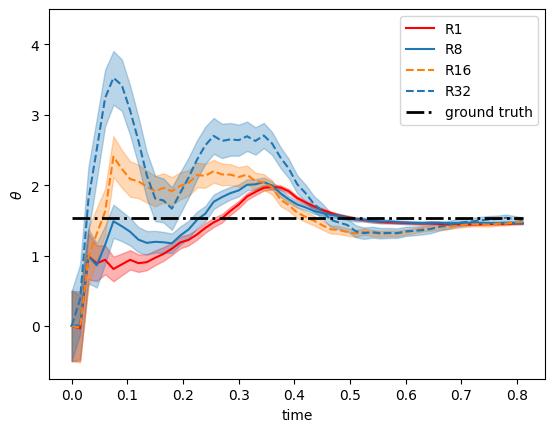} }
\subfloat[BART (Gaussian mask), $R_{d, 2}$]{
\includegraphics[trim=0 10 0 0, clip, width=0.3\textwidth]{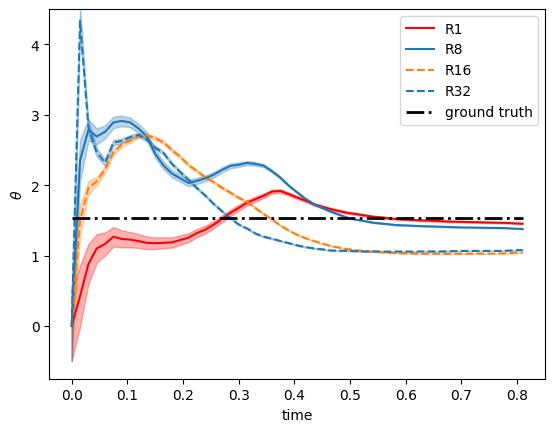} }\\
\subfloat[Spiral, $R_{d, 3}$]{
\includegraphics[trim=0 10 0 0, clip, width=0.3\textwidth]{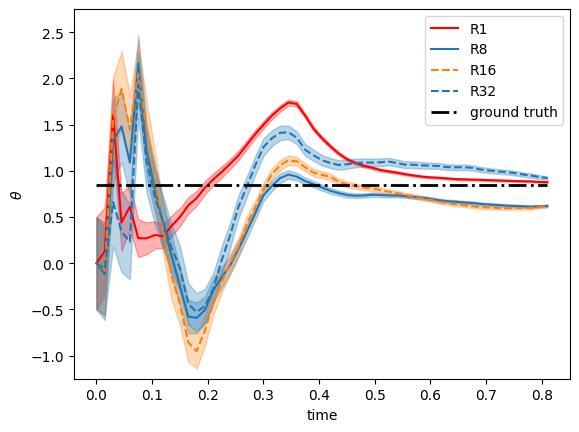} }
\subfloat[Gaussian, $R_{d, 3}$]{
\includegraphics[trim=0 10 0 0, clip, width=0.3\textwidth]{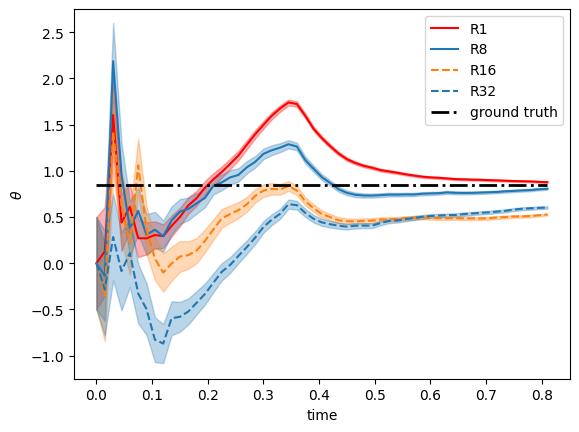} }
\subfloat[BART (Gaussian mask), $R_{d, 3}$]{
\includegraphics[trim=0 10 0 0, clip, width=0.3\textwidth]{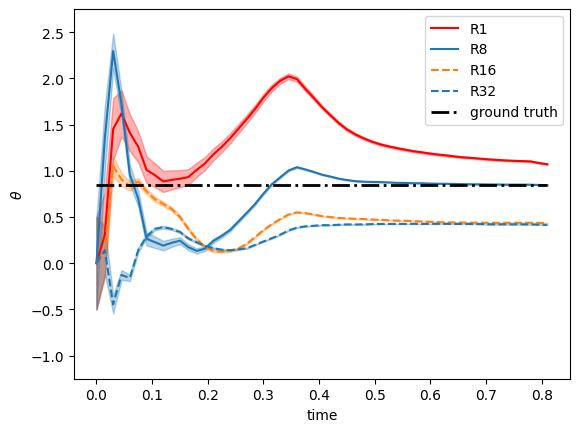} }\\

\caption{Evolution of the estimated parameters over the time of the simulation. The
dashed line indicates the true value.}
\label{fig:curves}
\end{figure}

\bigskip
\noindent\textit{Parameter curves of the Kalman filter.} We also provide examples of
the parameter curves of ROUKF in Figure \ref{fig:curves}, which show the development
of the estimated parameter over the length of the simulation. The shaded area
indicates the variance of the parameter. This parameter variance, for each of the
masks using our new method, is higher at the end of the simulation for higher $R$,
since less data is available. This effect is not observed for the BART-reconstructed
velocities, where the variance decreases similarly for all acceleration factors. Here
the variance is also much smaller than for the frequency-space based method,
indicating higher confidence in the estimated parameter. This higher confidence is
based on the fact that due to the reconstruction, the measurements contain more data
than in frequency space (as they now have as many voxels as the image space, as
compared to only the actually measured voxels). Since these data however originate
from the same subsampled measurements, they do not actually contain more information,
which leads to this misleading quantification of the uncertainty in the inverse problem.

Furthermore, it can be seen that for low $R$, the curves qualitatively follow the
shape of the curve for the fully sampled data, but with higher sampling rates, the
choice of the mask not only quantitatively but also qualitatively impacts the
evolution of the parameter estimation. The curves from our method also differ
considerably from those using BART, which develop erratic spikes in either direction
as the acceleration factor increases. In comparison to that, the curves for each mask
remain more similar to each other with increasing $R$.

It can also be seen that in all the cases, the value of the estimated parameter is
mostly constant after roughtly $t = 0.5$, and the parameter variance does not decrease
considerably after that. This indicates that for both frequency-space data and
reconstructed data, the relevant data is located entirely within systole.

In summary, for the synthetic data,
\begin{itemize}
\item using the inverse problem in frequency space significantly reduces the error
\item there are notable differences in the qualitative and quantitative estimation of
    the parameters between the different masks
\item the method is robust to different values of the $venc$ and using estimated
    magnitudes and background phases of the data.
\end{itemize}

\subsection{Phantom data}

\begin{figure}[htbp!]
\centering
\subfloat[Gaussian mask]{
\includegraphics[trim=30 10 0 30, clip, width=0.5\textwidth]{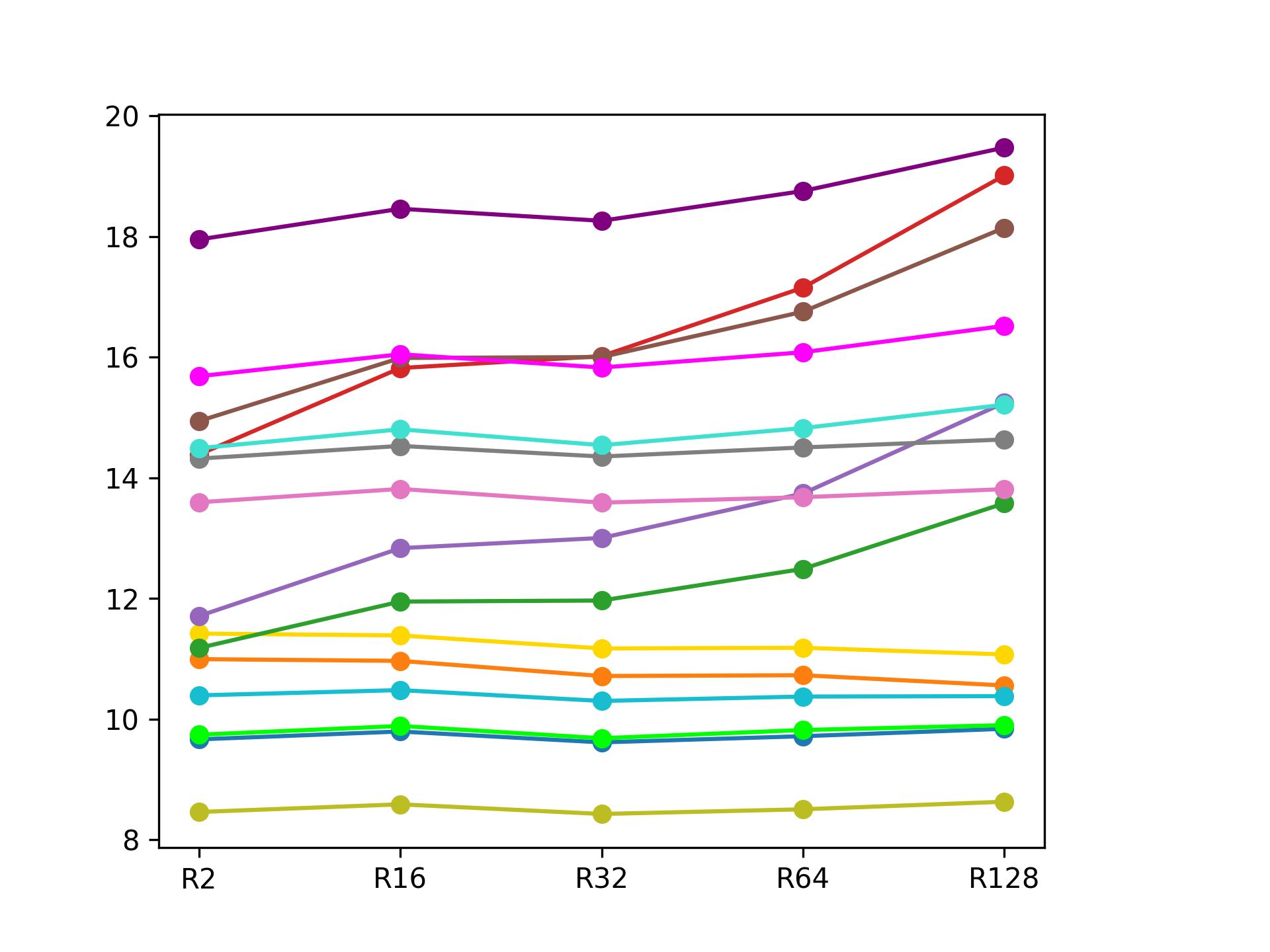}}
\subfloat[Spiral mask]{
\includegraphics[trim=30 10 0 30, clip, width=0.5\textwidth]{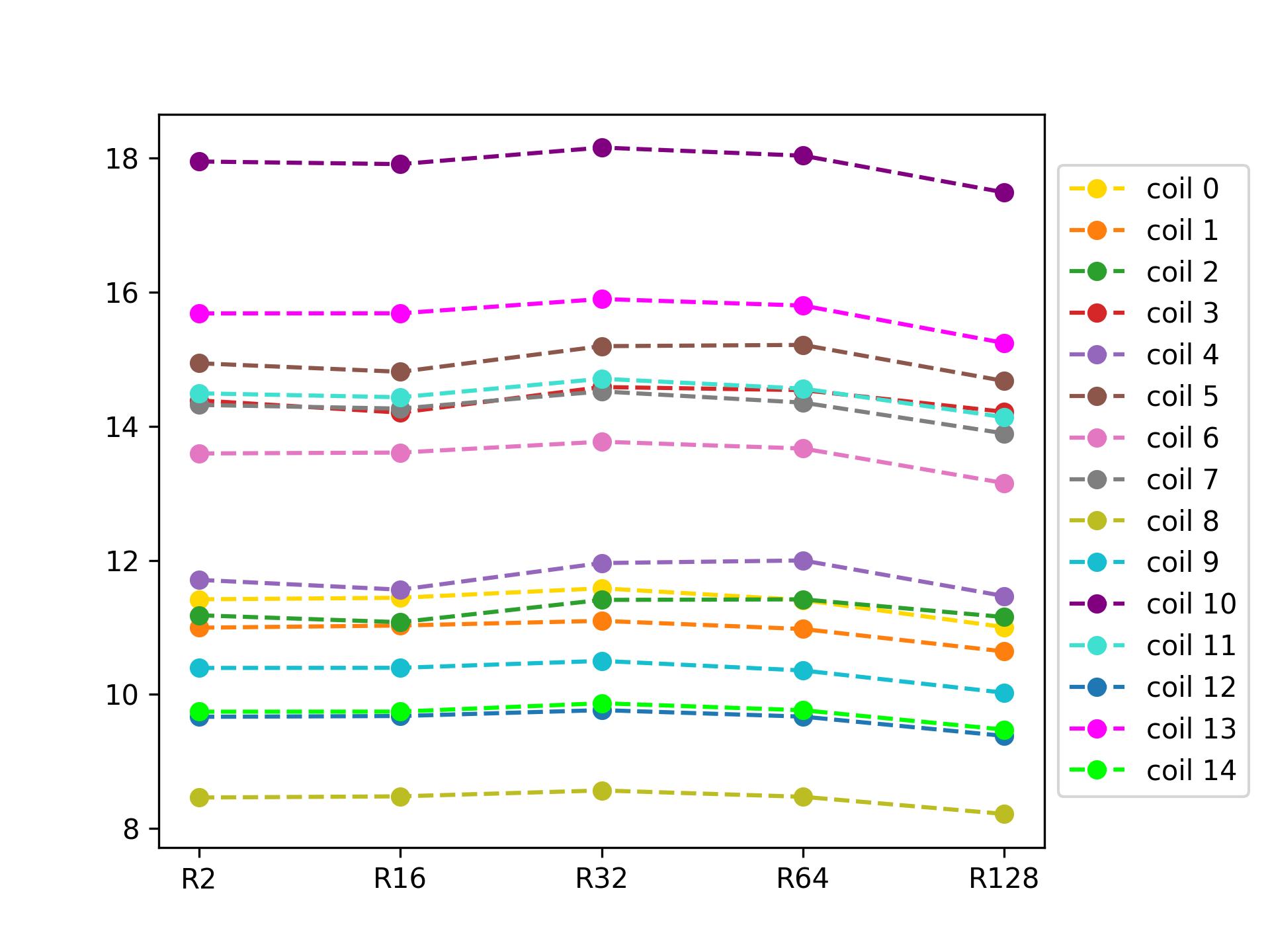}}
\caption{Estimation of the standard deviation of the signal noise per coil, for
various acceleration factors}
\label{fig:noise}
\end{figure}

First, we consider the quality of the estimation of the standard deviation of the
noise in the case of the phantom data. The results per mask and acceleration factor
are shown in Figure \ref{fig:noise}. While the estimation of the noise with the spiral
mask remains close to constant for the acceleration factors, the estimations with the
Gaussian mask show major variations for some of the coils. This is likely due to the
spatial variation of the signal strength for some of the coils and the fact that the
Gaussian mask does not uniformly cover the space. Nonetheless this approach provides a
good approximation of the standard deviation of the noise even for high acceleration factors.

For comparing the reconstructed flows, we are using the same error metric
\ref{eq:error} as in Section \ref{sec:synth_results}, with $\bs{u}_{ref}$ the solution
of a forward problem with the parameters estimated from $R = 2$.

\begin{figure}[!hbtp]
\centering

\includegraphics[trim=10 0 40 40, clip, width=0.6\textwidth]{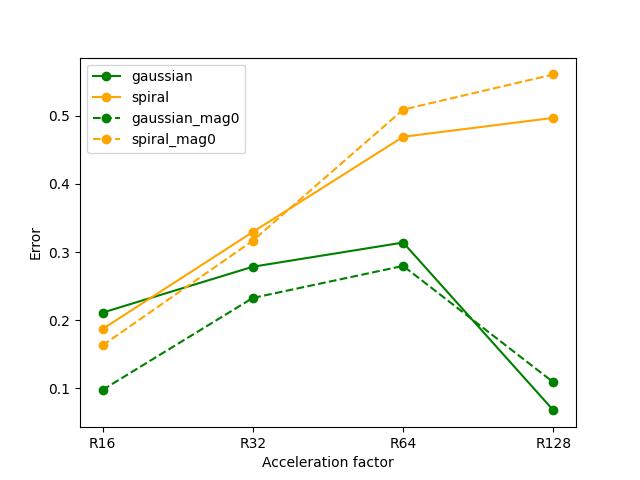}

\caption{Error values for different acceleration factors for both masks. The "mag0"
labelled results use the magnitude from the acquisition without a velocity encoding gradient.}
\label{fig:phantom_error}
\end{figure}

The error values can be seen in Figure \ref{fig:phantom_error}. A considerable
difference between the behaviour of the two masks is apparent. Both achieve similar
results for an acceleration factor of $R=16$, but the error increases higher for the
spiral mask than for the Gaussian mask. The Gaussian mask also shows a decrease in
error from $R=64$ to $R=128$, which may be a result of differences in the noise
estimation, or otherwise an outlier. Both choices for the magnitude in the observation
show similar results, implying robustness to using a magnitude from a different
velocity direction.

\begin{figure}[!hbtp]
\centering
\subfloat[Gaussian mask]{
    {
\includegraphics[trim=10 80 0 0, clip, width=\textwidth]{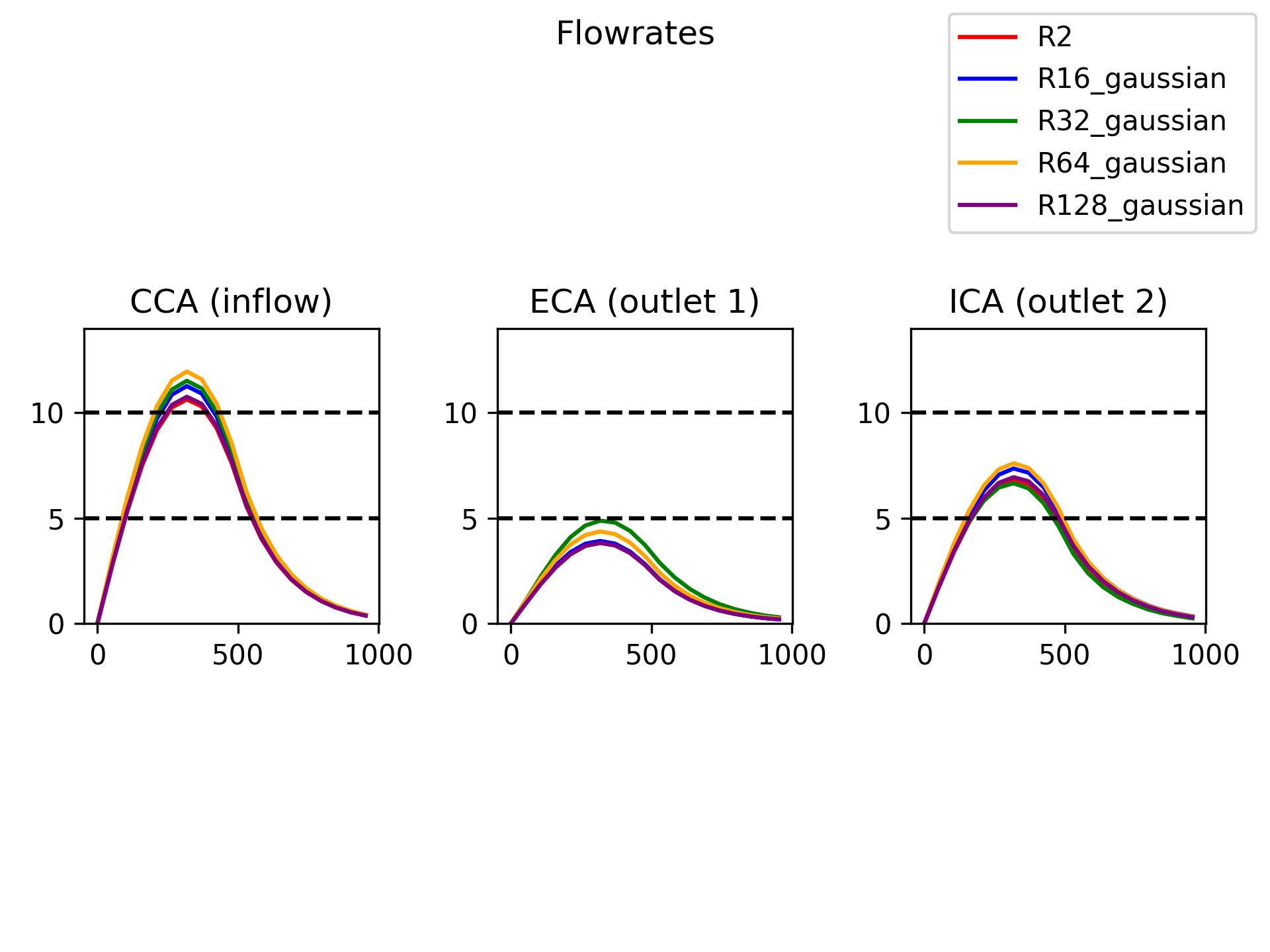} } }\\
\subfloat[Spiral mask]{
    \includegraphics[trim=10 80 0 100, clip, width=\textwidth]{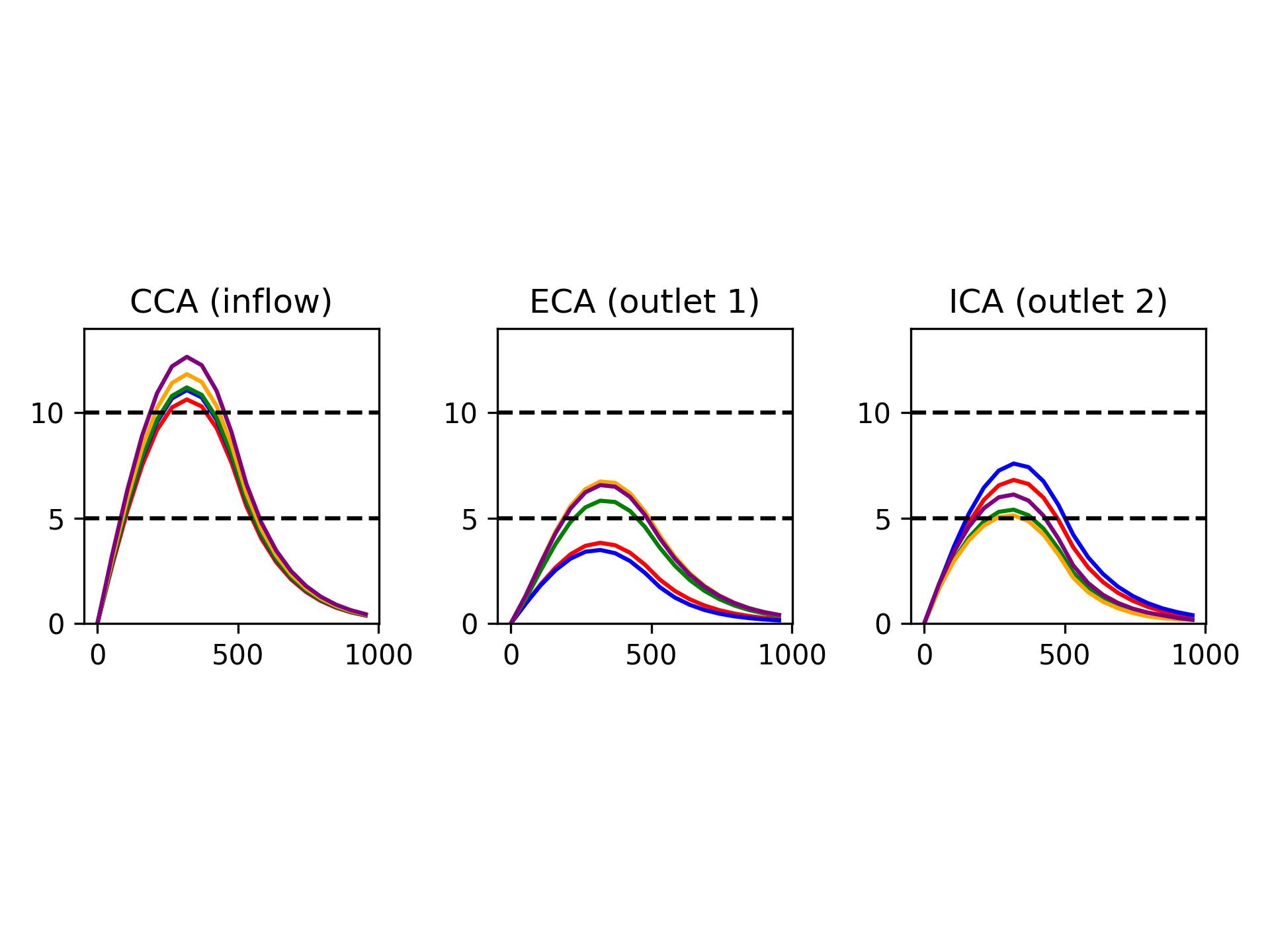}
}

\caption{Flow rates at the inlet and each of the outlets of the forward simulations
based on the estimated parameters.}
\label{fig:phantom_flow_rates}
\end{figure}

\begin{figure}[!hbtp]
\centering
\subfloat[$U$]{
    {
\includegraphics[trim=10 10 40 10, clip, width=0.47\textwidth]{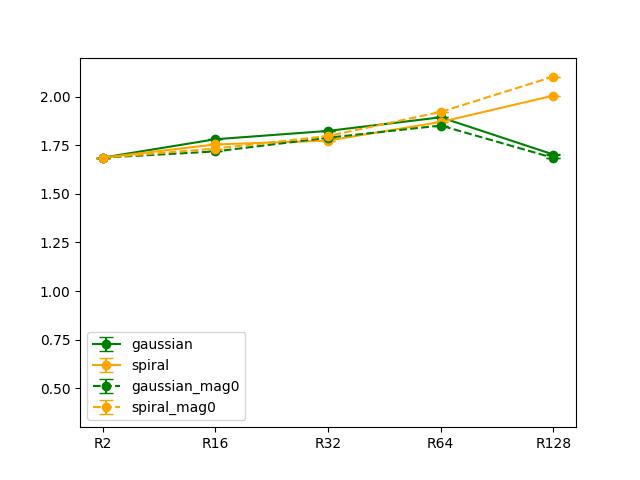} } }
\subfloat[$R_{p, 2}$]{
    \includegraphics[trim=10 10 40 10, clip, width=0.47\textwidth]{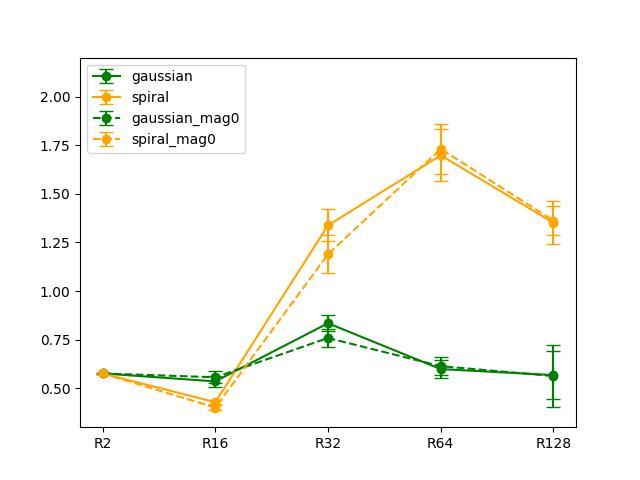}
}

\caption{Estimated parameter values for different $R$ for each mask. Errorbars depict
    the parameter variance, multiplied by a factor of 1000 for visibility. The "mag0"
labelled results use the magnitude from the acquisition without a velocity encoding gradient.}
\label{fig:phantom_params}
\end{figure}

In Figure \ref{fig:phantom_flow_rates}, we show the flow rates of the reconstructed
flows at the inlet and each of the outlets. It can be seen that both masks tend to
overestimate the inflow with increasing acceleration factor $R$. For the Gaussian
mask, the flow split between the two outlets remains fairly consistent as $R$
increases, whereas the spiral mask no longer accurately depicts the flow split from
$R=32$ upwards. This can also be observed by looking at the estimated parameter values
in Figure \ref{fig:phantom_params}. Both masks show a very consistent estimation of
the inflow amplitude $U$, with differences only for $R=128$, but the estimation of the
resistance boundary condition $R_{p, 2}$ shows larger changes for different
acceleration factors. The spiral mask especially overestimates $R_{p,2}$ for higher
acceleration factors, leading to a lack of distinction of the outlets in the flow rates.

\begin{figure}[!hbtp]
\centering
\subfloat[Gaussian mask]{
    {
\includegraphics[trim=5 5 15 5, clip, width=0.7\textwidth]{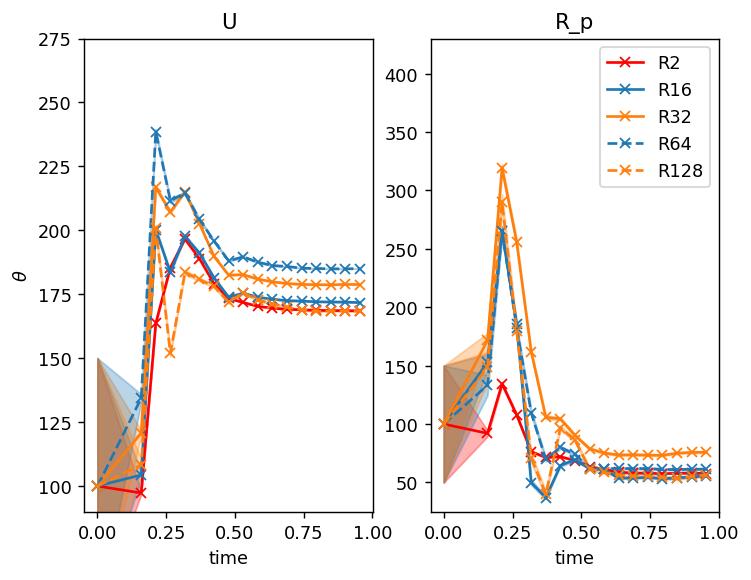} } }\\
\subfloat[Spiral mask]{
    \includegraphics[trim=5 5 15 5, clip, width=0.7\textwidth]{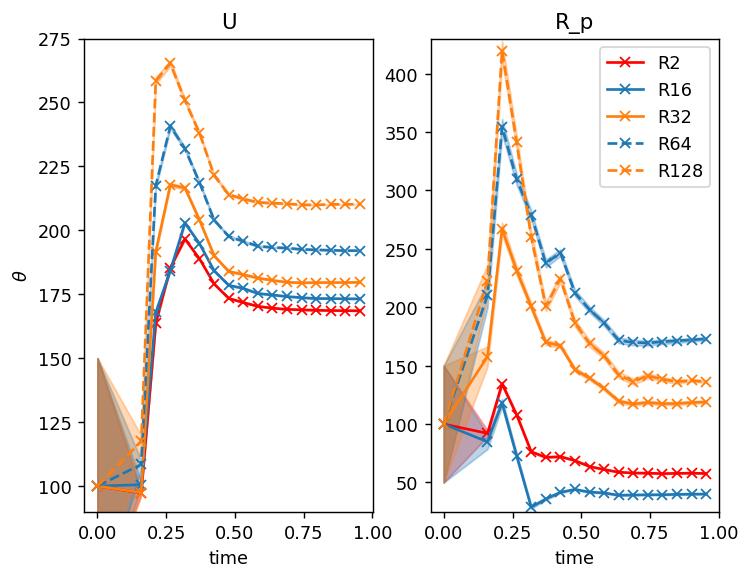}
}

\caption{Evolution of the estimated parameters over the time of the simulation.
Crosses mark times where corrections according to the measurements are made.}
\label{fig:phantom_curves}
\end{figure}

Examples of the ROUKF parameter curves are provided in Figure
\ref{fig:phantom_curves}. It can be seen that the curves for the Gaussian and spiral
masks have a similar shape with some qualitative differences. The acceleration factor
seems to have a higher effect on the parameter evolution curve for the spiral mask
than the Gaussian mask. It can also be seen that the parameter variance decreases
drastically with the first measurement, with only a very small variance afterward
despite continuing changes in the parameter values.

In summary, for the phantom data,
\begin{itemize}
\item the method successfully estimates the parameters with time-dependent and
    estimated magnitudes, background phases, and noise level
\item using the gaussian mask, the estimation is reasonably accurate even for very
    high acceleration factors (R=128)
\item the results are quantitatively and qualitatively different for the different
    masks, with the Spiral mask failing to distinguish between the two outlets for
    higher acceleration factors.
\end{itemize}

\section{Discussion}
The numerical results are consistent between the synthetic dataset and the phantom
dataset. In both cases, our approach robustly estimates boundary conditions of the
flow even with high undersampling rates and low $venc$ values and outperforms the
conventional approach of reconstructing velocity measurements from undersampled data
through Compressed Sensing.

The results also show significant differences between the two studied sampling
patterns, with qualitatively and quantitatively different outcomes of the inverse
problem. The sampling patterns appear to differ in their estimation of different
parameters, indicating that the choice of sampling pattern may depend on which
parameter is of higher interest.

In related research, in \cite{muller_reduced-order_2018}, Fourier-transformed
time-varying measurements (for only a single spatial point) to estimate parameters of
a one-dimensional fluid flow model, but this approach (neither the model nor the
inverse problem) does not utilize the potential of using velocity data of the entire
space in the frequency domain. Additionally, \cite{zhao_direct_2016} used an Unscented
Kalman Filter with undersampled frequency-space data to retrieve T2 mappings for MRI
brain scans. However, they are relying on an assumption that the phase is both known
and constant, hence only utilizing the magnitude of the k-space images, as well as
estimating the state rather than a number of parameters. Additionally, they only
experimented with relatively small acceleration factors (up to a factor of 8).

A limitation of the method is that it requires an accurate estimate of the magnitude
and the background phase of the magnetization. However, since the magnitude can be
scanned without motion encoding gradients, a highly sampled image of the magnitude is
quick to acquire. We have also shown that it was possible -- at least in the phantom
dataset -- to acquire the background image highly sampled for approximations of the
magnitude and background phase, while undersampling the remaining three directions. As
we have shown some robustness of the method with regards to flaws in the magnitude,
one could also acquire the magnitude only at diastole, and use this for all time steps
to cut down on acquisition time, or correct it using the model itself. Additionally,
the noise present in the reconstructed magnitude and background phase is currently not
considered by the Kalman filter. This will be subject of future research.

Another limitation of this paper is that subject data is not included, which would
pose additional complexities such as the respiratory motion of the vessels due to the
breathing of the patient and the potential large size of patient-specific meshes. We
plan to address this aspect in future research.

\section{Conclusion}
We have proposed a new formulation for the inverse problem for parameter estimation in
fluid flow problems using undersampled frequency-space MRI data and demonstrated it
using a Reduced-Order Unscented Kalman Filter. This method outperforms the results of
parameter estimation using velocity data reconstructed with Compressed Sensing,
especially for high undersampling rates and different $venc$ values.

The choice of subsampling mask is shown to have a strong influence on the result of
the estimation of some parameters. Future work could therefore address the challenge
of finding optimal sampling masks for certain parameters.

\section*{Acknowledgments}
C.B. and M.L. acknowledge the funding from the European Research Council (ERC) under
the European Union's Horizon 2020 research and innovation program (grant agreement No
852544 - CardioZoom).

\appendix
\section{Numerical solution method of the forward problem}
Here we detail the algorithm used to solve the incompressible Navier-Stokes equation
with Windkessel boundary conditions for the forward problem.

\begin{algorithm}[!hbtp]
\caption{Fractional step algorithm with a modified semi-implicit Windkessel model coupling}
\label{alg:semi}
Given the initial conditions $\bs{u}^0 = \bs u(0) \in V_{\Gamma_{w},h}$ and
$\pi^0_1,\dots,\pi_N^0 \in \mathbb{R}$,  perform for $j>0$, with $t^{j}  = j \tau$: \\

\textbf{1. Viscous Step:} Find the tentative velocity $\tilde{\bs{u}}^{n} \in
V_{\Gamma_{w},h}$ such that:
\begin{equation}
    \begin{cases}
        \displaystyle \tilde{\bs{u}}^{j}|_{\Gamma_{in}} = \bs{u}_{inlet}(t^{j}) \\
        \displaystyle \frac{\tau}{\rho} (\tilde{\bs{u}}^{j} , \bs{v} )_{\Omega_h} +
        \rho ( \bs{u}^{j-1} \cdot \nabla \tilde{\bs{u}}^{j} , \bs{v}  )_{\Omega_h} +
        \frac{\rho}{2} ( (\nabla \cdot \bs{u}^{j-1}) \tilde{\bs{u}}^{j}, \bs{v}  )_{\Omega_h}
        + (\delta \bs{u}^{j-1} \cdot \nabla \tilde{\bs{u}}^j , \bs{u}^{j-1} \cdot
        \nabla \bs{v}  )_{\Omega_h} \\
        \displaystyle + 2\mu (\epsilon(\tilde{\bs{u}}^{j}) ,
        \epsilon(\bs{v}))_{\Omega_h} + \sum_{\ell=1}^K \frac{\rho}{2} |
        \bs{u}^{j-1}\cdot \mathbf{n} |_{-} (\tilde{\bs{u}}^j,\bs{v})_{\Gamma_\ell} =
        \frac{\tau}{\rho} (\bs{u}^{j-1} , \bs{v})_{\Omega_h}
    \end{cases}
\end{equation} \label{eq:viscousstep}
for all $\mathbf{v} \in V_{\Gamma_{in}\cup \Gamma_{w} ,h}$, and $|x|_{-}$ denoting the
negative part of $x$.

\textbf{2. Projection-Windkessel Step:} Compute $\tilde{Q}^{j} = \int_{\Gamma_\ell}
\tilde{\bs{u}}^{j} \cdot \bs{j}$. Find $p^{j} \in Q_h$ such that:
\begin{equation} \label{eq:pressproj}
    \displaystyle \frac{\tau}{\rho} (\nabla p^{j} , \nabla q)_{\Omega_h} +
    \sum_{\ell=1}^K \frac{ \overline{p^{j}}_{\Gamma_\ell}
    \ \overline{q}_{\Gamma_\ell}}{\gamma_\ell} +  \epsilon
    \sum_{\ell=1}^K(\mathcal{T}(\nabla p^j), \mathcal{T}(\nabla q))_{\Gamma_\ell} =
    \sum_{\ell=1}^K \bigg( \tilde{Q}^{j} + \frac{\alpha_\ell
    \pi_\ell^{j-1}}{\gamma_\ell}  \bigg) \overline{q}_{\Gamma_\ell} - (\nabla \cdot
    \tilde{\bs{u}}^{j} , q)_{\Omega_h} ,
\end{equation}
for all $q \in Q_h$ and with $\overline{(\cdot)}_{\Gamma_{\ell}} =
\frac{1}{Area(\Gamma_\ell)} \int_{\Gamma_\ell} (\cdot) ds$ and
$\mathcal{T}(\mathbf{f}) = \mathbf{f} - (\mathbf{f}\cdot\mathbf{j}) \mathbf{j}$.

\textbf{3. Velocity correction Step:} Find $\bs{u}^{j} \in [L^2(\Omega_h)]^3 $ such that:
$$\displaystyle (\bs{u}^{j},\mathbf{v})_{\Omega_h} = (\tilde{\bs{u}}^{j} -
\frac{\tau}{\rho} \nabla p^{j},\mathbf{v})_{\Omega_h}
$$
for all $\mathbf{v} \in [L^2(\Omega_h)]^3$

\textbf{4. Update-Windkessel Step:} Set $P^{j}_\ell = \overline{p^{j}}_{\gamma_\ell} $
and compute $\pi_\ell^{j} \in \mathbb{R}$ as:
$$ \pi_\ell^{j} = \big ( \alpha_\ell - \frac{\alpha_\ell \beta_\ell}{\gamma_\ell} \big
) \ \pi_\ell^{j-1} + \frac{\beta_\ell}{\gamma_\ell} \  P^{j}_\ell \ , \ \ell = 1,...,K $$
\end{algorithm}
\section{ROUKF algorithm}
Here we detail the ROUKF algorithm adapted from \cite{moireau-chapelle-11}. Let us
first consider the notation $[\bs{Z}^{(*)}]$ as the matrix with the column-wise
collection of vectors $\bs{Z}_{(1)},\bs{Z}_{(2)},\dots$.

Define  the  \textit{canonical sigma-points}  $\bs{I}_{(1)},\dots,\bs{I}_{(2p)} \in
\mathbb{R}^p$ such that
\begin{equation*}
\bs{I}_{(i)} =
\begin{cases}
    \sqrt{p}\bs{e}_i, & \text{for } 1\leq i \leq p\\
    -\sqrt{p}\bs{e}_{i-p}, & \text{for } p+1 \leq i \leq 2p
\end{cases}
\end{equation*}
where the vectors $\bs{e}_i$ form the canonical base of of $\mathbb{R}^p$.
Moreover, define the weight $\alpha = \frac{1}{2p}$.

We denote by $\bs{\hat X}^n_-,\bs{\hat X}^n_+ \in\mathbb{R}^r$ a priori (model prediction) and a
posteriori (corrected by observations) estimates of the true state
$\bs{X}^n\in\mathbb{R}^r$. In the
semi-implicit coupled
3D-0D fractional step Algorithm~\ref{alg:semi}, the state consists
in the velocity field $\bs u^n$ and the Windkessel pressures $\pi^n_\ell$.
Estimates of all unknown parameters are summarized by the corresponding a
priori and a posteriori vectors $\bs{\hat\theta}^n_-,\bs{\hat\theta}^n_+\in\mathbb{R}^p$.
The discretized forward model is written as $\bs{X}^n =
\mathcal{A}^n(\bs{X}^{n-1},\bs{\theta}^{n-1})$,
$\mathcal{A}^n$ denoting the model operator.

For given values of the initial condition $\bs{\hat{X}}^0_+ = \bs{X}^0 \in \mathbb{R}^r$,
the initial expected value of the parameters $\bs{\hat{\theta}}^0_+ =
\bs{\theta}^0 \in\mathbb{R}^p$ and its covariance matrix $\bs{P}^0$, perform
\begin{itemize}
    \begin{subequations}\label{eq:ParamSEIK}
    \item \textbf{Initialization:}  initialize the sensitivities as
        \begin{multline}
            \bs{L}_{\theta}^{0} = \sqrt{\bs{P}^0} \ \text{ (Cholesky factor)},
            \bs{L}_X^{0} = \bs{0} \in \mathbb{R}^{r \times
            p},  \bs{U}_{0} =\bs{P}_\alpha\equiv
            \alpha[\bs{I}^{(*)}][\bs{I}^{(*)}]^T\label{eq:algo-estim-1}
        \end{multline}

        \noindent\hspace{-1.2cm} Then, for $n=1, \cdots, N_T$:
    \item \textbf{Sampling:} generate $2p$ particles from the current state and
        parameter estimates, i.e.~for $i=1,\ldots,2p$:
        \begin{multline}
            \shoveright{\hspace{0.1cm}
                \begin{cases}
                    \bs{\hat{X}}^{n-1}_{(i)} = \bs{\hat{X}}^{n-1}_+ +
                    \bs{L}_X^{n-1}(\bs{C}^{n-1})^{T}\bs{I}_{(i)}, \\
                    \bs{\hat{\theta}}^{n-1}_{(i)} = \bs{\hat{\theta}}^{n-1}_+ +
                    \bs{L}^{n-1}_{\theta}(\bs{C}^{n-1})^{T}\bs{I}_{(i)}
            \end{cases}}\label{eq:algo-estim-2}
        \end{multline}
        with $\bs{C}^{n-1}$ the Cholesky factor of $(\bs{U}^{n-1})^{-1}$.
    \item \textbf{Prediction:} propagate each particle with the forward model
        and compute an a priori state prediction:
        \begin{multline}
            \shoveright{\hspace{0.1cm}
                \begin{cases}
                    \bs{\hat{X}}^n_{(i)} =
                    \mathcal{A}^{n}(\bs{\hat{X}}^{n-1}_{(i)},\bs{\hat{\theta}}^{n-1}_{(i)}),\quad
                    \bs{\hat{\theta}}^{n}_{(i)} =\bs{\hat{\theta}}^{n-1}_{(i)},  \quad
                    i=1,\dots,2p \\
                    \bs{\hat{X}}^n_- = E_\alpha([(\bs{\hat{X}}^{n})^{(*)}]) \equiv
                    \alpha \sum_{i=1}^{2p} \bs{\hat{X}}^n_{(i)}  \\
                    \bs{\hat{\theta}}^n_- = E_\alpha([(\bs{\hat{\theta}}^n)^{(*)}])
            \end{cases}}\label{eq:algo-estim-3}
        \end{multline}
    \item \textbf{Correction:}  compute a posteriori estimates based on
        measurements for state and parameters, using
        the $i$-th
        particle innovation
        $\bs{\Gamma}^n_{(i)}$:
        \begin{multline}
            \shoveright{\hspace{0.1cm}
                \begin{cases}
                    \bs{L}^n_X = \alpha [(\bs{\hat{X}}^n)^{(*)}]  [\bs{I}^{(*)}]^{T}
                    \\
                    \bs{L}_\theta^n = \alpha [(\bs{\hat{\theta}}^n)^{(*)}]
                    [\bs{I}^{(*)}]^{T} \\ 
                    \bs{L}_\Gamma^n =  \alpha [(\bs{\Gamma}^n)^{(*)}]  [\bs{I}^{(*)}]^{T} \\
                    \bs{U}^n =  \bs{P}_\alpha +  (\bs{L}^n_\Gamma)^{T} \bs{W}^{-1}
                    \bs{L}^n_\Gamma\ ,  \quad \bs{P}_\alpha =  \alpha [\bs{I}^{(*)}]
                    [\bs{I}^{(*)}]^{T}\\
                    \bs{\hat{X}}^n_+ = \bs{\hat{X}}^n_- - \bs{L}^n_X
                    (\bs{U}^n)^{-1}(\bs{L}^n_\Gamma)^{T} \bs{W}^{-1}
                    E_\alpha([(\bs{\Gamma}^n)^{(*)}]) \\
                    \bs{\hat{\theta}}_+^n = \bs{\hat{\theta}}^n_- - \bs{L}^n_\theta
                    (\bs{U}^n)^{-1}(\bs{L}^n_\Gamma)^{T}
                    \bs{W}^{-1}E_\alpha([(\bs{\Gamma}^n)^{(*)}]) \\
                    \bs{P}_{\theta}^{n} = \bs{L}_{\theta}^{n} (\bs{U}^n)^{-1}
                    (\bs{L}_{\theta}^{n})^T\\
                    \bs{P}_{X}^{n} = \bs{L}_{X}^{n} (\bs{U}^n)^{-1}(\bs{L}_X^{n})^T
            \end{cases}}\label{eq:algo-estim-4}
        \end{multline}
    \end{subequations}
\end{itemize}
with $\bs{W}=\sigma_y^2 \mathbb{I}$.







\section*{Ethics Statement}
None

\bibliographystyle{elsarticle-num}
\bibliography{biblio_merged.bib,biblio_flow.bib}

\end{document}